\documentclass{amsart}
\usepackage{amssymb,amscd}
\newcommand{\nc}{\newcommand}
\nc{\Prnongen}{{\textsl{Pr}}_{\kappa,\mathrm{nondeg}}}
\nc{\Center}{{\mathcal{Z}}}
\nc{\Wcat}{\BGG(\W_{\kappa}(\sg))}
\nc{\fmap}{\rightsquigarrow }
\nc{\fdomain}{{\mathcal{C}}}

\newcommand{\HW}{{\mathcal{HW}}}

\newcommand{\why}{{\bf L}}

\newcommand{\M}{{\bf{M}}}

\newcommand{\U}{{\mathcal{U}}}

\newcommand{\Dg}{{\bf D}}

\newcommand{\W}{{\mathcal{W}}}

\newcommand{\sI}{\bar{I}}
\newcommand{\1}{{\mathbf{1}}}
\newcommand{\Ln}[1]{L\bar{\mathfrak{n}}_{#1}}

\newcommand{\w}{{\textsl{y}}}
\newcommand{\sQ}{\bar{Q}}
\newcommand{\eW}{\widetilde{W}}

\newcommand{\gen}{{\rm{gen}}}
\newcommand{\sroots}{\bar{\Delta}}
\newcommand{\roots}{\Delta}
\newcommand{\proots}{\Delta_+}

\newcommand{\rroots}{\Delta^{\rm re}}
\newcommand{\iroots}{\Delta^{\rm im}}
\newcommand{\prroots}{\Delta_+^{\rm re}}

\newcommand{\nrroots}{\Delta_-^{\rm re}}
\newcommand{\F}{{\mathcal{F}}}

\newcommand{\sPi}{\bar{\Pi}}

\newcommand{\Cl}{{\mathcal{C}l}}
\newcommand{\Tg}{L^{\g}}
\newcommand{\Tf}{L^f}

\newcommand{\TW}{L}

\newcommand{\semiinf}{\frac{\infty}{2}}

\newcommand{\Ttot}{L^{{\rm tot}}}

\newcommand{\C}{{\mathbb C}}
\newcommand{\Z}{{\mathbb Z}}

\newcommand{\inv}{^{-1}}
\newcommand{\dual}[1]{{#1}^*}
\newcommand{\lam}{\lambda}
\newcommand{\Lam}{\Lambda}

\renewcommand{\*}{{\otimes}}
\newcommand{\+}{\mathop{\oplus}}
\newcommand{\h}{ {\mathfrak h}}
\newcommand{\g}{ {\mathfrak g}}

\renewcommand{\sb}{\bar{\mathfrak{b}}}
\newcommand{\che}{^{\vee}}

\newcommand{\bra}{{\langle}}
\newcommand{\ket}{{\rangle}}

\DeclareMathOperator{\nondeg}{nondeg}
\DeclareMathOperator{\tp}{top}
\DeclareMathOperator{\vac}{vac}

\DeclareMathOperator{\coker}{coker}
\DeclareMathOperator{\Aut}{Aut}

\DeclareMathOperator{\End}{End}

\DeclareMathOperator{\Hom}{Hom}

\DeclareMathOperator{\id}{id}
\DeclareMathOperator{\ad}{ad}
\DeclareMathOperator{\haru}{span}

\DeclareMathOperator{\Ext}{Ext}

\DeclareMathOperator{\ch}{ch}

\DeclareMathOperator{\rank}{rank}
\DeclareMathOperator{\gr}{gr}
\DeclareMathOperator{\height}{ht}

\newcommand{\sg}{  \bar{\mathfrak g}}
\newcommand{\sh}{\bar{ \h}}
\newcommand{\sn}{\bar{\mathfrak{n}}}
\newcommand{\snp}{\bar{\mathfrak{n}}_+}
\newcommand{\snn}{\bar{\mathfrak{n}}_-}

\newcommand{\sP}{\bar P}
\newcommand{\sW}{\bar{W}}
\newcommand{\srho}{ \bar{\rho}}

\newcommand{\sproots}{\bar{\Delta}_+}
\newcommand{\snroots}{\bar{\Delta}_-}

\renewcommand{\check}{^{\vee}}

\DeclareMathOperator{\im}{Im}

\newcommand{\A}{{\mathcal A}}

\newcommand{\ud}[2]{{\genfrac{}{}{0pt}{}{#1}{#2}}}

\newcommand{\n}{{{\mathfrak{n}}}}

\newcommand{\BGG}{{\mathcal O}}

\newcommand{\Obj}{Obj}

\title[  Irreducible representations of
$\W$-algebras]{quantized
reductions
and
irreducible representations  of
$\W$-algebras}
\author{Tomoyuki Arakawa}
\address{Graduate school of Mathematics,  Nagoya
University,
Chikusa-ku, Nagoya, 464-8602, JAPAN}
\email{tarakawa@math.nagoya-u.ac.jp}
\dedicatory{Dedicated to Akihiro Tsuchiya on the occasion of his 60th
birthday}
\begin{document}
\theoremstyle{plain}
\newtheorem{Th}{Theorem}[subsection]
\newtheorem{MainTh}{Main Theorem}
\newtheorem{Pro}[Th]{Proposition}
\newtheorem{Lem}[Th]{Lemma}
\newtheorem{Co}[Th]{Corollary}

\newtheorem{Facts}[Th]{Facts}

\theoremstyle{definition}

\newtheorem{dfandpr}[Th]{Definition and Proposition}
\theoremstyle{remark}
\newtheorem{Def}[Th]{Definition}
\newtheorem{Rem}[Th]{Remark}
\newtheorem{Conj}{Conjecture}
\newtheorem{Claim}{Claim}
\newtheorem{Notation}{Notation}
\newtheorem{Ex}[Th]{Example}

\newcommand{\st}{{\mathrm{st}}}
\maketitle
\begin{abstract}
We study the representations
of the $\W$-algebra
$\W(\sg)$
associated to an arbitrary  finite-dimensional 
simple Lie algebra
$\sg$
via
the 
quantized Drinfeld-Sokolov  reductions.
The characters of
irreducible representations
of $\W(\sg)$ with
non-degenerate highest weights are expressed by
Kazhdan-Lusztig polynomials.
The irreduciblity conjecture 
of Frenkel, Kac and Wakimoto
is proved
completely for the ^^ ^^ $-$" reduction
and partially for  the ^^ ^^ $+$" reduction.
In particular,
the existence of the minimal series representations
($=$ the modular invariant representations)
of $\W(\sg)$ is proved.
\end{abstract}
\section{Introduction}
Since introduced by Zamalodchikov \cite{Z},
the symmetry by $\W$-algebras has been 
significantly important in conformal field theories (\cite{BS1, BS2}).
However,
not much  is known about the representation theory of $\W$-algebras.
In this paper
we 
study the representations of the  $\W$-algebra $\W(\sg)$
associated to an arbitrary simple Lie algebra $\sg$,
and
determine the characters of its irreducible highest weight representations
under certain conditions.
Those characters are expressed by Kazhdan-Lusztig polynomials.
As a consequence,
we prove the conjecture of
Frenkel, Kac and Wakimoto (\cite{FKW}) on
the existence of the minimal series representations
($=$ the modular invariant representations)
of $\W(\sg)$.

\smallskip

Let $\sg$
be a finite-dimensional 
simple Lie algebra
with a triangular decomposition $\sg=\sn_-\+ \sh\+ \sn_+$.
Let $\g=\sg\* \C[t,t\inv]
\+\C K\+ \C \Dg$  be the affine Lie algebra associated
to  $\bar{\mathfrak{g}}$.
Let  
$\W_{\kappa}(\sg)$
be the $\W$-algebra
associated to 
$\sg$ at level $\kappa-h\che$,
defined by Feigin and Frenkel
via the quantized Drinfeld-Sokolov reduction \cite{FF_W,FB}.
We have:
\begin{align}\label{eq:Zhu-iso-intro}
 \A(\W_{\kappa}(\sg))
\cong \Center(\sg),
\end{align}
where $\A(\W_{\kappa}(\sg))$
is the Zhu algebra of $\W_{\kappa}(\sg)$ (\cite{FZ})
and $\Center(\sg)$ is the center of the universal enveloping algebra
 $U(\sg)$,
see Theorem \ref{Th:Zhu-Center}.
By \eqref{eq:Zhu-iso-intro},
the irreducible highest weight representations
of $\W_{\kappa}(\sg)$
are parameterized by central characters 
of $\Center(\sg)$. Let
 $\gamma_{\bar{\lam}}:
\Center(\sg)\rightarrow \C$
be the central character
 defined as the evaluation at the Verma module of
$\sg$ of highest weight $\bar{\lam}\in \dual{\sh}$.
Let $\why(\gamma_{\bar{\lam}})$
be the irreducible representation of $\W_{\kappa}(\sg)$
of highest weight $\gamma_{\bar{\lam}}$,
that is,
the irreducible 
representation $\W_{\kappa}(\sg)$ corresponding to 
$\gamma_{\bar{\lam}}$.
To determine the  characters 
of all 
$\why(\gamma_{\bar{\lam}})$ is,
certainly, 
one of the
most  important problems in representation theory of $\W$-algebras.

\smallskip

Let
$\BGG_{\kappa}$ be
the Bernstein-Gelfand-Gelfand category
of $\g$ of level $\kappa-h\che$.
Let $H^{\bullet}_{+}(V)$,
$V\in \BGG_{\kappa}$,
be the cohomology associated to the
quantized Drinfeld-Sokolov reduction
defined by Feigin and Frenkel (\cite{FF_W}). 
Thus,
$H^{\bullet}_+(V)=H^{\semiinf +\bullet}(\Ln{+},
V\* \C_{\chi_+})$,
where
$\Ln{+}=\sn_+\* \C[t,t\inv]\subset \g$
and $\C_{\chi_+}$
is a certain one-dimensional representation of $\Ln{+}$.
Then,
the correspondence $V\fmap H^{i}_{+}(V)$
$(i\in \Z)$
defines a family of functors from $\BGG_{\kappa}$
to the category of $\W_{\kappa}(\sg)$-modules
(\cite{FF_W}).
However,
because
the description of
this functor is rather complicated in general,
Frenkel, Kac and Wakimoto \cite{FKW} introduced another functor
$V\fmap H^{i}_{-}(V)$
$(i\in \Z)$
 from $\BGG_{\kappa}$
to the category of $\W_{\kappa}(\sg)$-modules,
defined by
$H^{\bullet}_-(V)=H^{\semiinf +\bullet}(\Ln{-},
V\* \C_{\chi_-})$.
Here,
$\Ln{-}=\sn_-\* \C[t,t\inv]\subset \g$
and $\C_{\chi_-}$
is again a certain one-dimensional representation of $\Ln{-}$.
It turns out that
the corresponding functor 
\begin{align*}
 \text{$V\fmap H^{0}_{-}(V)$}
\end{align*}
indeed has nice properties:
it is exact,
it sends Verma modules to Verma modules,
and simple modules to simple modules
under certain  conditions.

\smallskip

We now describe our
results more precisely.
Let $\h$
be the Cartan subalgebra of $\g$
and let $\dual{\h}\ni \Lam \mapsto \bar{\Lam}
\in \dual{\sh}$ be the restriction.
Let $M(\Lam)$
be the Verma module of $\g$
of highest weight $\Lam\in \dual{\h}$,
$L(\Lam)$
its unique simple quotient.
Let $\sproots $
be the set of positive roots of $\sg$,
identified with the subset of the set $\prroots$
of positive real roots
of $\g$
in the standard way.
\begin{MainTh}[{Theorem \ref{Th:irr-} }]\label{MainTh1}
 Let $\Lam$ be non-degenerate 
$($that is,
$\bra \Lam,\bar{\alpha}\che\ket\not\in \Z$
for all $\bar{\alpha}\in \sproots)$
and non-critical $($that is, $\kappa=
\bra \Lam+\rho,K\ket \ne 0)$.
Then,
$H_-^0(L(\Lam))\cong \why(\gamma_{\bar{\Lam}})$.
\end{MainTh}
Note that the condition that $\Lam$ is non-degenerate 
does not imply that
the integral Weyl group $W^{\Lam}$ of $\Lam$
is finite
(and indeed it is infinite when $\Lam$ is admissible,
see \cite{KW2}).
We also show that
$H_-^0(M(\lam))\cong \M(\gamma_{\bar{\lam}})$,
where $\M(\gamma_{\bar{\lam}})$
is the Verma module of $\W_{\kappa}(\sg)$
of highest weight $\gamma_{\bar{\lam}}$,
see Theorem \ref{Th:image_of_Verma_Module}.
Thus,
combined with our previous result \cite{A},
Main Theorem \ref{MainTh1}  
gives the  character 
of $\why(\gamma_{\bar{\lam}})$ with any
 non-degenerate  highest weight $\gamma_{\bar{\lam}}$,
see Theorem \ref{Th:ch-formula}.
Here,
an infinitesimal character $\gamma_{\bar{\lam}}$ ($\bar{\lam}\in 
\dual{\sh}$) is called non-degenerate if $\bra
\bar{\lam},\bar{\alpha}\che\ket \not\in \Z$
for all $\bar{\alpha}\in \sproots$.

\smallskip

By applying  Main Theorem \ref{MainTh1}  to 
non-degenerate admissible weights (\cite{KW2}),
the irreduciblity conjecture of Frenkel, Kac and Wakimoto
\cite{FKW} is proved (for the ``$-$''-reduction).
Therefore,
combined with our previous result \cite{A},
the existence of modular invariant representations
of $\W(\sg)$ is proved.

\medskip
Now
the ``$+$''-reduction
is not as simple as the ``$-$''-reduction.
However,
we have the following theorem.
\begin{MainTh}[{Theorem \ref{Th:simple+}}]\label{MainTH2}
 Suppose that $\Lam\in \dual{\h}$ is non-critical and
satisfies the following condition:
\begin{align*}
 \bra \Lam,\alpha\che\ket \not \in \Z\quad
\text{for all $\alpha\in \{-\bar{\alpha}+n\delta;
\bar{\alpha}\in \sproots, 1\leq n\leq \height \alpha\}$}.
\end{align*}
Then,
$H_+^0(L(\Lam))\cong \why(\gamma_{\overline{t_{-\srho\che}\circ \Lam}
})$.
\end{MainTh}

\smallskip

This paper is organized as follows.
In Section 2,
we collect the necessary information about the affine Lie algebra
$\g$ and the  BRST complexes
associated to quantized Drinfeld-Sokolov reductions.
In Section 3, we recall the definition of $\W$-algebra
$\W_{\kappa}(\sg)$
and the collect necessary information about its structure.
In Section \ref{section:Zhu},
we prove that the Zhu algebra $\A(\W_{\kappa}(\sg))$ of 
$\W_{\kappa}(\sg)$
is canonically isomorphic to
the center $\Center(\sg)$ of 
$U(\sg)$
(Theorem \ref{Th:Zhu-Center}).
In Section 5,
we define Verma modules $\M(\gamma_{\bar{\lam}})$ 
of $\W_{\kappa}(\sg)$ of highest weight $\gamma_{\bar{\lam}}$ and
its simple quotient $\why(\gamma_{\bar{\lam}})$.
The duality structure of modules over $\W$-algebras
is also discussed.
Finally,
in Section 6,
we prove  Main Theorems.

\bigskip

{\em Acknowledgments. }
I would like to thank Edward Frenkel for 
valuable  discussions.
\section{Preliminaries
on the affine Lie algebra
and the BRST complex}
In this section we collect the necessary information about
the  affine Lie algebra
and the BRST complexes associated to quantized Drinfeld-Sokolov
reductions.
Our basic references in this section
are the textbooks
\cite{KacBook,MPBook,FB}.
\subsection{The affine Lie algebra $\g$}
In the sequel,
we fix 
a simple
finite-dimensional complex Lie algebra $\sg$ and a
Cartan subalgebra $\sh\subset \sg$. Let $\sroots$ denote the set
of roots,
$\sPi$  a basis of $\sroots$,
$\sproots$  the set of positive roots,
and $\snroots=-\sproots$.
This gives the triangular
decomposition $\sg=\snn\+\sh\+\snp$.
Let $\sQ$ denote the root lattice,
$\sP$ the weight lattice,
$\sQ\che$ the coroot lattice
and
$\sP\che$ the coweight lattice.
Let $\srho$ be the half  sum of
positive roots,
$\srho\che$ the half  sum of positive coroots.
For $\alpha\in \sproots$,
the number
$\bra \alpha,\srho\che\ket$
is called the
{\em height} of
$\alpha$ and denote by
$\height \alpha$.
Let $\sW$ be the Weyl group of $\sg$,
$w_0$ the longest element of $\sW$.

Let $(~,~)$ be the normalized invariant inner product
of $\sg$.
Thus,
$(~,~)=\frac{1}{2h\check}$Killing form,
where
$h\check$ is the dual Coxeter number of $\sg$.
We identify $\sh$ and $\dual{\sh}$ using
the form.
Then,
$\alpha\che
=2\alpha/(\alpha,\alpha)$,
$\alpha\in \sroots$.

Let $l=\rank \sg$,
$\sI=\{1,2,\dots,l\}$.
Choose a basis $\{J_a;a\in
\sI\sqcup
\sroots\}$  of $\sg$
such that
$J_{\alpha}\in \sg_{\alpha}$,
$(J_{\alpha},J_{-\alpha})=1$
and $(J_{\alpha})^t=J_{-\alpha}$
($\alpha\in \sroots$).
Here, $\sg\ni X\mapsto X^t \in \sg$
is the Chevalley anti-automorphism.
Let $c_{a,b}^c$ be the
structure constant with respect to
this basis;
$[J_a,J_b]=\sum_{c}c_{a,b}^c J_c$.
Then,
$c_{\alpha,\beta}^{\gamma}
=-c_{-\alpha,-\beta}^{-\gamma}$
($\alpha,\beta,\gamma\in \sproots$).

\smallskip

Let
$\g=\sg\*\C[t,t^{-1}]\+\C K \+ \C \Dg$
be the affine Lie algebra associated to
($\sg$,$(~,~)$),
where
$K$ is its central element and $\Dg$ is the degree
operator
(\cite{KacBook}).
The bilinear form $(~,~)$
is naturally extended from $\sg$ to $\g$.
Set $X(n)=X\* t^n$,
$X\in \sg$,
$n\in \Z$.
The subalgebra $\sg\* \C\subset \g$
is naturally identified with $\sg$.

Fix
the  triangular decomposition
$\g=\g_-\+\h\+\g_+$
in the standard way.
Thus,
\begin{align*}
&\h=\sh\+\C K\+\C \Dg,\\
&\text{$\g_-=\snn\* \C[t\inv]\+ \sh\* \C[t\inv]t\inv
\+ \snp\*\C[t\inv]t\inv$,}\\
&\text{$\g_+=\snn\* \C[t]t\+ \sh\* \C[t]t
\+ \snp\*\C[t]$.}
\end{align*}
Let
$\dual{\h}=\dual{\sh}\+\C \Lam_0\+\C \delta$
be the  dual of $\h$.
Here,
$\Lam_0$ and
$\delta$ are dual elements of $K$ and $\Dg$
respectively. For $\lam\in\dual{\h}$,
the number
$\bra \lam,K\ket$ is called the
{\em level of $\lam$}.
Let
$\dual{\h}_{\kappa}$
denote the
set of the
weights of level
$\kappa-h\che$:
\begin{align}
\dual{\h}_{\kappa}
=\{\lam\in \dual{\h};\bra \lam+\rho,K\ket=\kappa\},
\end{align}
where,
$\rho=\bar{\rho}+h\che
\Lam_0\in\dual{\h}$.
Let
$\bar{\lam}$ be the restriction of $\lam\in \dual{\h}$ to $\dual{\sh}$.

\smallskip

Let $\roots$ be the set of roots of $\g$,
$\roots_+$ the set of positive roots,
$\roots_-=-\roots_+$.
Then,
$\roots=\rroots\sqcup \iroots$,
where
$\rroots$
is
the set of real roots
and
$\iroots$ is the set of imaginary roots.
Let
$\Pi$ be the standard basis of $\rroots$,
$\rroots_{\pm}=\rroots\cap \roots_{\pm}$,
$\iroots_{\pm}=\iroots\cap \roots_{\pm}$.
Let $Q$ be the root lattice,
$Q_+=\sum_{\alpha\in \proots}\Z_{\geq 0} \alpha\subset
Q$.

Let
$W\subset GL(\dual{\h})$ be the Weyl group of $\g$
generated by the reflections $s_{\alpha}$,
$\alpha\in \rroots$,
defined by
$s_{\alpha}(\lam)=\lam-\bra \lam,\alpha\che\ket\alpha$.
Then,
$W=\sW\ltimes \sQ\che$.
Let $\eW=\sW\ltimes \sP\che$,
the extended Weyl group of $\g$.
For $\mu\in \sP\che$,
we denote the corresponding element
of $\eW$ by $t_{\mu}$.
Then,
\begin{align*}
t_{\mu}(\lam)=\lam+\bra\lam,K\ket\mu-
\left(\bra\lam,\mu\ket+\frac{1}{2}
|\mu|^2\bra\lam,K\ket
\right)\delta \quad (\lam\in \dual{\h}).
\end{align*}
Let $\eW_+=\{w\in \eW;\prroots\cap
w\inv(\nrroots)=\emptyset  \}$.
We have:
$\eW=\eW_+\ltimes W $.
The dot action of
$\eW$
on $\dual{\h}$ is defined by $w\circ \lam= w(\lam+\rho)-\rho$
($\lam\in \dual{\h}$).
We have the natural homomorphism
$\eW\rightarrow \Aut (\g)$
such that $w(\g_{\alpha})\subset \g_{w(\alpha)}$
($w\in \eW$, $\alpha\in \roots$).
\smallskip

For $\Lam\in \dual{\h}$,
let
\begin{align}
 \text{$R^{\Lam}=\{
\alpha\in \rroots;
\bra \Lam+\rho,\alpha\che\ket\in \Z\}$,
}
\end{align}$R^{\Lam}_+=R^{\Lam}\cap
\prroots$.
It is known that
$R^{\Lam}$ is a subroot
system of $\rroots$
(\cite{MPBook,KT1}).
Let 
\begin{align}
 \text{$W^{\Lam}
=\bra s_{\alpha};\alpha\in R^{\Lam}\ket\subset W$.
}
\end{align}
The Coxeter group $W^{\Lam}$
is called the {\em{integral Weyl group}} of $\Lam$.
We have:
\begin{align}
\text{$
R^{w\circ \Lam}=R^{\Lam}$
for all $w\in W^{\Lam}$.}
\end{align}

\subsection{The BGG category
of $\g$}
For a $\g$-module $V$
(or for  simply an $\h$-module $V$),
let $V^{\lam}=\{v\in V;
hv=\lam(h)v \text{ for }h\in \h\}$
be the weight space of weight $\lam$.
Let $P(V)=\{\lam\in \dual{\h};V^{\lam}\not=
\{0\}\}$.
If
$\dim V^{\lam}<\infty$ for all $\lam$,
then
we set
\begin{align}\label{eq:graded-dual}
V^*=\bigoplus_{\lam}\Hom_{\C}(V^{\lam},\C)
\subset \Hom_{\C}(V,\C).
\end{align}

\smallskip

Let
$\BGG_{\kappa}$ be the
full subcategory of the category of left $\g$-modules
consisting of objects $V$
such that
(1) $V$ is locally finite
over $\g_+$,
(2) $V=\bigoplus\limits_{\lam\in
\dual{\h}_{\kappa}}V^{\lam}$ and
$\dim_{\C}V^{\lam}<\infty$ for all $\lam$, (3)
there exists
a finite subset
$\{\mu_1,\dots,\mu_n\}\subset
\dual{\h}_{\kappa}$
such that
$P(V)\subset  \bigcup\limits_i \mu_i-Q_+$.

The correspondence
$V\rightsquigarrow V^*$
defines
the duality functor
in $\BGG_{\kappa}$.
Here, $\g$ acts on $V^*$
by
$(Xf)(v)=f(X^tv)$,
where
$X\mapsto X^t$ is  the Chevalley
antiautomorphism of $\g$.

Let $M(\lam)\in \BGG_{\kappa}$,
$\lam\in \dual{\h}_{\kappa}$,
be the Verma module of
highest weight $\lam$
and
$L(\lam)$ its unique simple
quotient.
Let $\BGG_{\kappa}^{[\Lam]}$,
$\Lam\in \dual{\h}_{\kappa}$,
be the
full subcategory of
$\BGG_{\kappa}$
whose objects have all their
local composition factors
isomorphic to $L(w\circ\Lam)$,
$w\in W^{\Lam}$.
By \cite{Kumer},
$\BGG_{\kappa}$ splits into the orthogonal
direct sum
$
\BGG_{\kappa}=\bigoplus\limits_{\Lam
\in \dual{\h}_{\kappa}/\sim}\BGG_{\kappa}^{[\Lam]}
$,
where
$\sim$ is the equivalent relation
defined by
$\lam\sim \mu\Leftrightarrow
\mu \in W^{\lam}\circ \lam$.
Orthogonal here means that
$\Ext^i_{\BGG_{\kappa}}(M,N)=0 $
for $M\in \BGG_{\kappa}^{[\Lam]}$,
$N\in \BGG_{\kappa}^{[\Lam']}$,
$i\geq 0$,
when $\Lam\ne \Lam'$ in $\dual{\h}_{\kappa}/\sim$.
\subsection{The Clifford algebra and the  Fock space}
Let 
\begin{align}
\text{$\Ln{\pm}
=\sn_{\pm}\*\C[t,t\inv]\subset \g$.
}\end{align}
We identify $\Ln{\pm}$ with $(\Ln{\mp})^*$
using the invariant bilinear form $(~,~)$
of $\g$.
Let
$\Cl$ be  the {\em Clifford algebra}
associated to
$\Ln{+}\+\Ln{-}
$
and its  symmetric bilinear form
defined by the identification $(\Ln{\pm})^*=\Ln{\mp}$.
Denote by
$\psi_{\alpha}(n)$,
$\alpha\in \sroots$,
$n\in \Z$,
the
generators of $\Ln{+}\+\Ln{-}
\subset \Cl$
which correspond to
the elements
$J_{\alpha}(n)$.
Then,
\begin{align*}
\{\psi_{\alpha}(m),{\psi}_{\beta}(n)\}=
\delta_{\alpha+\beta,0}\delta_{m+n,0}\quad
(\alpha,\beta\in \sroots, m,n\in \Z).
\end{align*}
Here,
$\{x,y\}=xy+yx$.
The algebra $\Cl$   contains
the Grassmann algebra
$\Lam(\Ln{\pm})
$ of $\Ln{\pm}$ as its subalgebra.
We have
$\Cl=\Lam(\Ln{+})\*\Lam(\Ln{-})$ as  $\C$-vector
spaces.
The action of $\eW$ on $\g$ naturally extends to $\Cl$:
in particular
$\sP\che$
acts as
\begin{align}
 \text{$t_{\mu}(\psi_{\alpha}(n))=\psi_{\alpha}(n-\bra \alpha,\mu\ket)$
\quad 
($\mu \in \sP\che$,
$\alpha\in \sroots$,
$n\in \Z$).
}
\end{align}

Let $\bar{\Cl}$
be the subalgebra of $\Cl$
generated by $\psi_{\alpha}(0)$,
$\alpha\in \sroots$.
Then,
$\bar{\Cl}$ is identified the
 Clifford algebra associated to 
the space $\sn_+\+ \sn_-$
and its natural non-degenerate bilinear form.
We have:
\begin{align}
 \bar{\Cl}\cong \Lam (\sn_+)\* \Lam (\sn_-)
\end{align}
as $\C$-vector spaces.

Let $\F(\Ln{\pm})$  be the irreducible
representation of $\Cl$ generated by a vector $\1$
such that
\begin{align}
 \text{$\psi_{\alpha}(n)\1=0$
\quad 
($\alpha \in \sroots$,
$n\in \Z$,
$\alpha+n\delta\in \prroots$).
}
\end{align}
Thus,
$\F(\Ln{\pm})=\Lam(\Ln{\mp}\cap \g_-)\*\Lam(\Ln{\pm}\cap \g_-)$
as  $\C$-vector spaces.
Let
\begin{align}
\F^p(\Ln{\pm})=
\sum\limits_{i-j=p}\Lam^i(\Ln{\mp}\cap \g_-)\*\Lam^j(\Ln{\pm}\cap \g_-)
\subset \F(\Ln{\pm})
\quad (p\in \Z).\label{eq:identification-of-F}
\end{align}
Then,
$\F(\Ln{\pm})=\sum_{p\in \Z}\F^p(\Ln{\pm})$.
\subsection{The BRST complex}
For $V\in \BGG_{\kappa}$,
let
$$C(\Ln{\pm},V)=
V\* \F(\Ln{\pm})=\sum_{i\in \Z}
C^i(\Ln{\pm},V),
\quad \text{where $C^i(\Ln{\pm},V)=V\*
\F^i(\Ln{\pm})$}.$$

Let
$d_{\pm}=d_{\pm}^{\st}+\chi_{\pm}\in \End C(\Ln{\pm},V)$,
where
\begin{align}
&d_{\pm}^{\st}=\sum_{\alpha\in
\sproots,n\in \Z}
J_{\pm \alpha}(-n){\psi}_{\mp \alpha}(n)
-\frac{1}{2}\sum_{\ud{\alpha,\beta,\gamma\in
\sproots}{ k+l+m=0}}
c_{\pm \alpha,\pm\beta}^{\pm\gamma}
{\psi}_{\mp \alpha}(k){\psi}_{\mp \beta}(l)
\psi_{\pm \gamma}(m),\\
&\chi_+=\sum_{\alpha\in
\sPi}{\psi}_{-\alpha}(1),\\
&\chi_-=\sum_{\alpha\in \sPi}\psi_{\alpha}(0).
\end{align}
Then,
$(d_{\pm}^{\st})^2=\chi_{\pm}^2=0$,
$\{d_{\pm}^{\st},\chi_{\pm}\}=0$
and
$d_{\pm} C^i(\Ln{\pm},V)\subset
C^{i+1}(\Ln{\pm},V)$.
In particular,
$d_{\pm}^2=0$.
For $V\in \Obj\BGG_{\kappa}$,
define
\begin{align}
H_{\pm}^{\bullet}(V)=H^{\bullet}
(C(\Ln{\pm},V),d_{\pm}).
\end{align}
It is called the
{\em cohomology 
of the BRST complex of the
quantized Drinfeld-Sokolov reduction
for $\Ln{\pm}$
associated to  $V$} (\cite{FF_W,FB,FKW}).
\section{$\W$-algebras}
In this section we collect necessary information
about $\W$-algebras.
The textbook \cite{FB}
and the paper \cite{FKW, FZ}
are our basic references
in this section.
\subsection{The vertex operator algebra
$C_{\kappa}(\sg)$}
Let 
\begin{align}
 \text{$V_{\kappa}(\g)=
U(\g)\*_{U(\sg\*\C[t]\+\C K\+
\C \Dg)}\C\in \Obj
\BGG_{\kappa}$.
}
\end{align}Here,
$\C$
is 
considered as a
$\sg\*\C[t]\+\C K\+
\C \Dg$-module
on which
$\sg\*\C[t]\+
\C \Dg$
acts trivially
and $K$ acts as $(\kappa-h\che)\id$.
It is called the
{\em universal affine vertex algebra
of level $\kappa-h\che$}
associated to $\sg$.
Let
\begin{align*}
C_{\kappa}(\g)
=C(\Ln{+},V_{\kappa}(\g))
=V_{\kappa}(\g) \* \F(\Ln{+}),
\end{align*}
and let $|0\ket=(1\*1)\* \1$ be its highest weight vector.
The space
$C_{\kappa}(\g)$
has 
a vertex  algebra structure,
see \cite[Chapter 14]{FB}.
Let $Y(v,z)\in
\End C_{\kappa}(\g)[[z,z\inv]]$
be the field corresponding to $v\in C_{\kappa}(\g)$.
Set
\begin{align*}
& X(z)=\sum_{n\in \Z}X(n)z^{-n-1}=Y(X(-1)|0 \ket ,z)\quad (X\in \sg)\\
&\psi_{\alpha}(z)=\sum_{n\in \Z}\psi_{\alpha}(n)z^{-n-1}
=Y(\psi_{\alpha}(-1)|0\ket,z)\quad (\alpha\in \sproots),\\
&
\psi_{-\alpha}(z)=\sum_{n\in \Z}\psi_{-\alpha}(n)z^{-n}=
Y(\psi_{-\alpha}(0)|0\ket,z)
\quad (\alpha\in \sproots).
\end{align*}
We have:
\begin{align}\label{eq:commutativiry_of_d_with_Y}
[d_+, Y(v,z)]=Y(d_+v,z)\quad
\text{for all $v\in C_{\kappa}(\g)$.
}
\end{align}

Let
\begin{align*}
L(z)=\sum_{n\in \Z}
L(n)z^{-n-2}=\Ttot(z)+\partial_z
\widehat{h}_{\srho\che}(z),
\end{align*}
where
\begin{align*}
&\Ttot(z)=
\sum_{n\in \Z}
\Ttot(n)z^{-n-2}=\Tg(z)+\Tf(z),\\
&\Tg(z)\text{ is the Sugawara filed of $\g$, }
\\&
\Tf(z)=-\sum_{\alpha\in\sproots}:\psi_{\alpha}(z)
\partial\psi_{-\alpha}(z):,\\
&\widehat{h}_{\srho\che}(z)=
\sum_{n\in \Z}\widehat{h}_{\srho\che}(n)z^{-n-1}
=
{\srho\che}(z)+\sum_{\alpha\in\sproots}
\height \alpha:\psi_{\alpha}(z)\psi_{-\alpha}(z):.
\end{align*}
Here,
$:~:$
is the normal ordering in the sence of \cite{FB}.
Then,
\begin{align}
&
[d_+, L(z)]=0,\label{eq:comm_Virasoro-element}\\
&
L(z)L(w)
\sim
\frac{c(\kappa)/2}{(z-w)^4}+
\frac{2}{(z-w)^2}L(w)+
\frac{1}{z-w}\partial L(w),
\label{eq:OPE-of-T}
\end{align}
where
\begin{align}\label{eq:cc}
c(\kappa)
&=l-12
\left(\kappa|\srho\che|^2
-2 \bra \srho,
\srho\che\ket+ \frac{1}{\kappa}|\srho|^2
\right).
\end{align}
The OPE
\eqref{eq:OPE-of-T} is equivalent to the
following commutation relations:
\begin{align*}
[L(m),L(n)]=(m-n)L(m+n)+\frac{
m^3-m}{12}\delta_{m+n,0}c(\kappa)\id.
\end{align*}
Note that
\begin{align}\label{eq:Visasoro}
\TW(-1)=\Ttot(-1),\quad
\TW(0)=\Ttot(0)-\widehat{h}_{\srho\che}(0),
\quad
\TW(1)=\Ttot(1)-2\widehat{h}_{\srho\che}(1).
\end{align}
In particular,
$\TW(0)$
acts on $C_{\kappa}(\g)$
semisimply  and
\begin{align}
Y(L(-1)v,z)=\partial_z Y(v,z).
\end{align}
Thus,
$C_{\kappa}(\g)$
has the  vertex operator algebra  structure  with Virasoro field $L(z)$.

\smallskip

Let
\begin{align*}
C_{\kappa}(\g)=\bigoplus_{\Delta\in
\Z}C_{\kappa}(\g)_{\Delta},
\end{align*}
be the eigenspace decomposition
with respect to the action of $\TW(0)$.
Thus,
\begin{align*}
C_{\kappa}(\g)_{\Delta}=\bigoplus_{\lam\atop
\bra \lam,\Dg+\srho\che\ket=-\Delta}
C_{\kappa}(\g)^{\lam}.
\end{align*}
Here,
$C_{\kappa}(\g)^{\lam}$ is the weight space of weight $\lam$
of $C_{\kappa}(\g)$ with respect to the natural action of $\h$.

For $v\in C_{\kappa}(\g)_{\Delta}$,
we expand
the corresponding field $Y(v,z)$
as
\begin{align}
Y(v,z)=\sum_{n\in \Z}Y_n(v)z^{-n-\Delta}
\label{eq:expansion_of_Y}
\end{align}
so $[L(0),Y_n(v)]=-n Y_n(v)$.
Thus,
for example,
$Y_n(v)=J_{\alpha}(n-\height \alpha)$,
$\psi_{\alpha}(n-\height \alpha)$,
$J_{-\alpha}(n+\height \alpha)$,
$\psi_{-\alpha}(n+\height \alpha)$
for $v=J_{\alpha}(-1)|0\ket$,
$\psi_{\alpha}(-1) |0\ket$,
$J_{-\alpha}(-1)|0\ket$,
$\psi_{-\alpha}(0)|0\ket$
respectively ($n\in \Z$,
$\alpha\in \sproots$).

\subsection{The
Zhu algebra of $C_{\kappa}(\g)$}
For a vertex operator algebra  $V$ in general,
one associates an  associative algebra
$\A(V)$
called the {\em Zhu algebra} (\cite{FZ}):
Let
\begin{align}
\text{$\U(V)=\bigoplus_{n\in \Z}\U(V)_n$,
\quad
$\U(V)_n=\{u\in \U(V);[L(0),u]=-nu \}$)
}
\end{align}
be the universal enveloping algebra
of $V$ in the sense of \cite{FZ}.
According to \cite{FZ},
one can define
$\A(V)$  as 
\begin{align}\label{eqdef:Zhu-algebra}
\A(V)=\U(V)_0/
\sum_{p>0}\U(V)_{-p}\U(V)_p,
\end{align}
(see \cite{NT} for a  proof of the above identification 
\eqref{eqdef:Zhu-algebra} of $\A(V)$).

The zero-mode algebra $\U_0(C_{\kappa}(\g))
\subset \U(C_{\kappa}(\sg))$
contains 
$t_{\srho\che}(U(\sg)\* \bar{\Cl})
\subset U(\g)\* \Cl$
generated by
\begin{align*}
 \text{$J_{\pm \alpha}(\mp \height
 \alpha)=t_{\srho\che}(J_{\pm\alpha}(0))$,
$\psi_{\pm \alpha}(\mp \height \alpha)=t_{\srho\che}(\psi_{\pm\alpha}(0))
$ ($\alpha\in \sproots$).}
\end{align*}
The following is clear by definition.
\begin{Pro}\label{Pro:Zhu_of_complex}
For any $\kappa\in \C$,
the Zhu algebra $\A(C_{\kappa}(\g))$
of $C_{\kappa}(\g)$
is canonically  isomorphic to
$t_{\srho\che}(U(\sg)\* \bar{\Cl})\cong 
U(\sg)\* \bar{\Cl}$.
\end{Pro}
\subsection{The tensor product decomposition of $C_{\kappa}(\g)$}
Following \cite{FZ},
we now recall the tensor product decomposition
of  $C_{\kappa}(\sg)$.
Set
\begin{align*}
\widehat{J}_{a}(z)=
\sum_{n\in \Z}\widehat{J}_{a}(n)z^{-n-1}
=J_a(z)+\sum_{\beta,\gamma\in \sproots}
c_{a,\beta}^{\gamma}:\psi_{\gamma}(z)\psi_{-\beta}(z):
\end{align*}
for $a\in \sI\sqcup \sroots$.
Let $C_{\kappa}(\g)_0$
be the subspace of
$C_{\kappa}(\g)$ spanned by the elements
\begin{align*}
\widehat{J}_{a_1}(-n_1)\dots
\widehat{J}_{a_p}(-n_p)
{\psi}_{-\alpha_{j_1}}
(-m_1)\dots {\psi}_{-\alpha_{j_q}}(-m_q)
|0\ket
\end{align*}
with
$a_i\in \sI\sqcup \snroots$,
$\alpha_{j_i}\in \sproots$,
$n_i, m_i\in \Z$.
Similarly,
let
 $C_{\kappa}(\g)'$
be the subspace of
$C_{\kappa}(\g)$ spanned by the elements
\begin{align*}
\widehat{J}_{\alpha_{i_1}}(-n_1)\dots
\widehat{J}_{\alpha_{i_p}}(-n_p)
{\psi}_{\alpha_{j_1}}
(-m_1)\dots {\psi}_{\alpha_{j_q}}(-m_q)
|0\ket
\end{align*}
with
$\alpha_{i_s},\alpha_{j_s}\in \sproots$,
$n_i, m_i\in \Z$.
It was shown in \cite{FB}, generalizing
the results of \cite{dBT},
that
\begin{align}
&\text{$C_{\kappa}(\g)'
$ and  $C_{\kappa}(\g)_0$
are vertex subalgebras of $C_{\kappa}(\g)
$,}
\label{eq:tensor-dec-1}
\\
& 
 d_+C_{\kappa}(\g)'\subset C_{\kappa}(\g)',\
d_+C_{\kappa}(\g)\subset C_{\kappa}(\g),\\
&\text{$C_{\kappa}(\g)=C_{\kappa}(\g)'\* C_{\kappa}(\g)_0$
as vertex algebras and complexes, }
\\
& H^i(C_{\kappa}(\g)')=
\begin{cases}
 \C &(i=0)\\
0&(i\ne 0),
\end{cases}\label{eq:vanishing-FB-1}\\
&
H^i(C_{\kappa}(\g)_0)=0\
(i\ne 0).
\label{eq:vanishing-FB-2}
\end{align}
In particular,
by K\"{u}nneth Theorem,
we have
\begin{align}\label{eq:iso-W-coho}
 H_+^i(V_{\kappa}(\g))=\begin{cases}
H^0(C_{\kappa}(\g)_0)
&(i=0)\\
0&(i\ne 0).
				    \end{cases}
\end{align}
\subsection{The Zhu algebra of $C_{\kappa}(\sg)'$ and
$C_{\kappa}(\sg)_0$}
The 
zero-mode algebra 
$\U_0(C_{\kappa}(\sg)')$
of $C_{\kappa}(\sg)'$
contains the algebra  $t_{\srho\che}(U(\sn_+)\* \Lam (\sn_+))
\cong U(\sn_+)\* \Lam (\sn_+)$
generated by
\begin{align*}
 \widehat{J}_{\alpha}(-\height \alpha),
\psi_{\alpha}(-\height \alpha)
\quad (\alpha\in \sproots).
\end{align*}
Here,
the vector space $U(\sn_+)\* \Lam (\sn_+)$
is considered as an algebra 
such that
(1) the natural maps
$U(\sn_+) \hookrightarrow
U(\sn_+)\* \Lam (\sn_+)$
and $\Lam (\sn_+)\hookrightarrow
U(\sn_+)\* \Lam (\sn_+)$ are algebra embeddings,
and (2) $[u,\omega]=\ad u(\omega)$
$(u\in U(\sn_+), \omega \in \Lam(\sn_+))$.

Similarly,
$\U_0(C_{\kappa}(\sg)_0)$
contains a subalgebra $t_{\srho\che}(U(\sb_-)\* \Lam (\sn_-))
\cong U(\sb_-)\* \Lam (\sn_-)$
generated by
\begin{align*}
\widehat{J}_i(0),
 \widehat{J}_{-\alpha}(\height \alpha),
\psi_{-\alpha}(\height \alpha)
\quad (i\in \sI,\alpha\in \sproots).
\end{align*}
Here,
$\sb_-=\sn_-\+\sh\subset \sg$
and $U(\sb_-)\* \Lam (\sn_-)$
is considered as an algebra similarly as above.

We have the following proposition.
\begin{Pro}
Let $\kappa$ be any complex number.

 \begin{enumerate}
  \item $\A(C_{\kappa}(\sg)')\cong  U(\sn_+)\* \Lam (\sn_+) $.
  \item $\A(C_{\kappa}(\sg)_0)\cong  U(\sb_-)\* \Lam (\sn_-) $.
 \end{enumerate}
\end{Pro}

\subsection{The tensor product decomposition of $\A(C_{\kappa}(\g))$}

Let the differential
$d_+$
act on $\g(C_{\kappa}(\g))$
by
$d_+J(v,f)=J(d_+v,f)$.
This makes
 $\U(\W_{\kappa}(\sg))$
and $\A(\W_{\kappa}(\sg))$
complexes.
Similarly,
the spaces
$\U(C_{\kappa}(\sg)')$,
$\U(C_{\kappa}(\sg)_0)$, $\A(C_{\kappa}(\sg)')$
and $\A(C_{\kappa}(\sg)_0)$
have the natural structure of complexes.

One can apply 
the argument of \cite[Section 14.2]{FB}
to prove the 
following
proposition.
\begin{Pro}\label{Pro:Zhu-easy}
Let $\kappa$ be any complex number.

 \begin{enumerate}
  \item $\A(C_{\kappa}(\g))
=\A(C_{\kappa}(\g)')\*
\A(C_{\kappa}(\g)_0)$
as complexes.
 \item
$H^i(\A(C_{\kappa}(\g)'))=
\begin{cases}
\C&(i=0)\\
0&(i\ne 0).
\end{cases}$
\item $H^i(\A(C_{\kappa}(\g)_0))=0$ $(i\ne 0)$.
 \end{enumerate}
\end{Pro}
By Proposition \ref{Pro:Zhu-easy},
we conclude as
\begin{align}
H^i(\A(C_{\kappa}(\g)))=
\begin{cases}
 H^0(\A(C_{\kappa}(\g)_0))&(i=0)\\
0&(i\ne 0).
\end{cases}
\label{eq:iso:Zhu-tensor-1}
\end{align}

\subsection{The $\W$-algebras}
Define
\begin{align}\label{eq:Def-of-W-algebras}
\W_{\kappa}(\sg)=H_+^0(V_{\kappa}(\g)).
\end{align}
By \eqref{eq:commutativiry_of_d_with_Y},
$Y$ descends to a map
\begin{align}\label{eq:field_of_W}
Y:\W_{\kappa}(\sg)\rightarrow
\End \W_{\kappa}(\sg) [[z,z\inv]].
\end{align}
Hence,
by \eqref{eq:comm_Virasoro-element},
$\W_{\kappa}(\sg)$
has a vertex operator algebra  structure 
with the Virasoro field 
$L(z)$
with the central charge $c(\kappa)$.
The vertex operator algebra  $\W_{\kappa}(\sg)$
is called the
{\em $\W$-algebra associated to $\sg$
at level $\kappa-h\che$}.
Note that, by \eqref{eq:iso-W-coho},
we have
\begin{align}\label{eq:realization_of_W_by_FB}
\W_{\kappa}(\sg)=H^0(C_{\kappa}(\g)_0)
\end{align}
as vertex algebras.

\subsection{A filtration of $\W_{\kappa}(\sg)$}
We now define a decreasing  filtration
\begin{align*}
\dots \supset G^p \W_{\kappa}(\sg)\supset 
G^{p+1}\W_{\kappa}(\sg)\supset \dots
\supset G^1\W_{\kappa}(\sg)=\{0\}
\end{align*}
of $\W_{\kappa}(\sg)$
such that 
\begin{align}\label{eq:multi-gr-W}
 G^p \W_{\kappa}(\sg)\cdot G^q \W_{\kappa}(\sg)
\subset G^{p+q}\W_{\kappa}(\sg)
\end{align}
and
the corresponding graded vertex algebra
$\gr^G \W_{\kappa}(\sg)=\bigoplus_p G^p
\W_{\kappa}(\sg)/G^{p+1}\W_{\kappa}(\sg)$
is  commutative.
Here,
the left-hand-side of \eqref{eq:multi-gr-W}
denotes the  span of the vectors
$Y_k(v_1)v_2$
($v_1\in G^p\W_{\kappa}(\sg)$,
$v_2\in G^q \W_{\kappa}(\sg)$,
$k\in \Z$).

\smallskip

Set
\begin{align}\label{eq:def-of-filtration}
G^{p}C_{\kappa}^n(\g)_0
=\bigoplus_{\lam\in
\dual{\h}
\atop
\bra \lam,
\srho\che\ket
\geq p-n}C_{\kappa}^n(\g)_0^{\lam}
\subset C_{\kappa}^n(\g)_0\quad (p\leq n+1),
\end{align}
where 
$C_{\kappa}^n(\g)_0=C_{\kappa}(\g)_0\cap
C_{\kappa}^n(\g)$. Then,
\begin{align}
&\dots \supset
G^{p}C_{\kappa}^n(\g)_0
\supset G^{p+1}C_{\kappa}^n(\g)_0
\supset \dots
\supset G^{n+1}C_{\kappa}^n(\g)_0=\{0\},\\
&C_{\kappa}^n(\g)_0=\bigcup_p G^{p}C_{\kappa}^n(\g)_0
\\
& d^{\st}_+G^{p}C_{\kappa}^n(\g)_0\subset
G^{p+1}C_{\kappa}^{n+1}(\sg)_0,\quad
\chi_+ G^{p}C_{\kappa}^{n}(\sg)_0
\subset
G^{p}C_{\kappa}^{n+1}(\sg)_0.
\end{align}
Let $E_r\Rightarrow H^{\bullet}(C_{\kappa}(\g)_0,d_+)$
be the corresponding spectral sequence.
By definition,
 $E_1=H^{\bullet}(C_{\kappa}(\g)_0,\chi_+)$.

Let $p_-=\sum_{i\in \sI}
\frac{(\alpha_i,\alpha_i)}{2}f_i$.
Then,
there exists a basis 
$\{P_i\}_{i\in \sI}$ of 
$\ker {\ad p_-}_{|\sn_-}
\subset \sn_-$
such that
$[{\bar{\rho}\che},P_i]=-d_i P_i$
($i\in \sI$),
where $d_i$ is the $i$-th exponent of $\sg$.
It is known that
$\ker {\ad p_-}_{|\sn_-}
=\haru \{P_i\}_{i=1}^l$
is a maximal abelian subalgebra of $\sg$.

It was shown in \cite[14.2.8]{FB}
that 
\begin{align}
&H^i(C_{\kappa}(\g)_0,\chi_+)=0\ (i\ne 0),
\\& 
\begin{array}{ccc}
\C[\widehat{P}_{1}(-n_1),
\dots, \widehat{P}_{l}(-n_l)]_{n_1,\dots,n_l\geq 1}
& \cong &H^0(C_{\kappa}(\g)_0,\chi),\\
f&\mapsto &f {|0\ket }
\end{array}
\label{eq;basis-of-gr-W}
\end{align}
where
 $\widehat{P}_i(n)$ is the linear combination
of $\widehat{J}_{\alpha}(n)$
corresponding to $P_i$.
In particular,
the 
spectral sequence
collapses at $E_1=E_{\infty}$.

Define
\begin{align}
 G^p \W_{\kappa}(\sg)
=\im: H^{0}(G^pC_{\kappa}(\g)_0)
\hookrightarrow H^{0}(C_{\kappa}(\g)_0)=\W_{\kappa}(\sg).
\end{align}
Then,
by definition,
\begin{align}
 \text{$\gr^G \W_{\kappa}(\sg)
=E_{\infty}=H^0(C_{\kappa}(\g)_0,\chi_+)$.}
\end{align}
Since
\begin{align}\label{eq:comp-action-fil}
G^p C^n_{\kappa}(\g)_0\cdot 
G^q C^m_{\kappa}(\g)_0 \subset G^{p+q} C^{n+m}_{\kappa}(\g)_0,
\end{align}
the graded vector space $\gr^G \W_{\kappa}(\sg)
$
carries a vertex algebra structure 
Here,
the left-hand-side of \eqref{eq:comp-action-fil}
denotes the  span of the vectors
$Y_k(v_1)v_2$
($v_1\in G^p C^n(\sg)_0$,
$v_2\in G^q C^m(\sg)_0$,
$k\in \Z$).
Since
$\gr^G \W_{\kappa}(\sg)\cong H^0(C_{\kappa}(\sg)_0,\chi_+)$
is as  vertex algebras,
$\gr^G \W_{\kappa}(\sg)$
is a commutative vertex algebra.
Clearly,
\begin{align}
 \A(\gr^G \W_{\kappa}(\sg))\cong
\A(H^0(C_{\kappa}(\g),\chi_+))= 
\C[ \widehat{P}_1(d_1), \widehat{P}_2(d_2), \dots,
\widehat{P}_l(d_l)
].\label{eq-gr-zhu-1}
\end{align}

Let
$W_i$ be the cocycle in $C_{\kappa}(\g)_0$
corresponding to $\widehat{P}_{i}(-1)|0\ket$.
Set
\begin{align*}
W_i(z)=\sum_{n\in \Z}W_i(n)z^{-n-d_i-1}
=Y(W_i,z).
\end{align*}
Then,
$W_i(n)|0\ket =0$ ($n\geq -d_i$),
$W_i=W_i(-d_i-1)|0\ket $
 ($i\in \sI$).
By \eqref{eq;basis-of-gr-W},
we see that 
the set
\begin{align*}
 \left\{
W_{i_1}(-d_{i_1}-n_1)\dots
W_{i_m}(-d_{i_m}-n_m)|0\ket;
\begin{array}
{l}
1\leq i_1\leq \dots \leq i_m\leq l,
n_p\geq 1,
\\n_p\geq n_{p+1}\text{ if
$i_p=i_{p+1}$ }
\end{array}\right\}.
\end{align*}
forms a $\C$-basis of $\W_{\kappa}(\sg)$.
Thus,
the vertex algebra $\W_{\kappa}(\sg)$ is freely generated
by the fields $W_1(z),\dots,W_{l}(z)$
in the sense of \cite{SK}.
\subsection{A filtration of $H^0(\A(C_{\kappa}(\g)_0))$}
Similarly as above,
we can define a filtration 
$\{G^p H^0(\A(C_{\kappa}(\g)_0))\}$
on the $\C$-algebra $H^0(\A(C_{\kappa}(\g)_0))$.
Let $\gr^G H^0(\A(C_{\kappa}(\g)_0))$
be
the corresponding graded algebra.
Then,
we have
\begin{equation}
 \begin{aligned}
 \gr^G H^0(\A(C_{\kappa}(\g)_0))&=
H^0(\A(C_{\kappa}(\g)_0), \chi_+)
\\&=\C[ \widehat{P}_1(d_1), \widehat{P}_2(d_2), \dots,
\widehat{P}_l(d_l)
].
\end{aligned}\label{eq-gr-zhu2}
\end{equation}
\subsection{The category $\Wcat$.}
Assume $\kappa\ne 0$.
For a $\U(\W_{\kappa}(\sg))$-module
$V$,
let
$V_{h}^{\gen}$ be the generalized eigenspace of
$L(0)$ of eigenvalue $h\in \C$: 
\begin{align}
\text{$V_{h}^{\gen}
=\{v\in V;
(L(0)-h)^nv=0 \text{ for $n\gg 0$}\}$.
}\end{align}
Let $\Wcat $
be the full subcategory of $\U(\W_{\kappa}(\sg))$-modules
consisting of objects  $V$ such that
\begin{align*}
(1) &\text{ there exists a finite set
$\{h_1,\dots ,h_r\}$ in $\C$ such that
 $V=\bigoplus\limits_{h\in \bigcup_i (h_i+\Z_{\geq 0})}V^{\gen}_{h}$,
}\\
(2)&\text{ $\dim V^{\gen}_{h}<\infty $
for all $h\in \C$.}
\end{align*}
It is easy to see that
the category $\Wcat$ is abelian.

For an object $V$ in $\Wcat$,
we define the normalized character 
$\ch V$
by
\begin{align}
 \ch V=\sum_{h\in \C}q^{h-\frac{c(\kappa)}{24}
}\dim_{\C}
V_{h}^{\gen}.
\end{align}

\subsection{The functors}\label{section:functors}
By the definition \eqref{eq:Def-of-W-algebras}
of $\W_{\kappa}(\sg)$,
the space $H_+^i(V)$,
$V\in \Obj \BGG_{\kappa}$
is a naturally $\W_{\kappa}(\sg)$-module.
It is also true that
$H_-^i(V)$,
$V\in \Obj \BGG_{\kappa}$,
has a $\W_{\kappa}(\sg)$-module structure.
The action of
$\U(\W_{\kappa}(\sg))$
on 
$H_-^{\bullet}(V)
$ is twisted by the element
${\w}=w_0t_{-\srho\che}=
t_{\srho\che}w_0\in \eW$:
\begin{align}
u\cdot v=\w(u)v\quad 
(u\in \U(\W_{\kappa}(\sg)),
v\in H_-^{\bullet}(V) ).
\end{align}
The above action is well-defined since
\begin{align}\label{eq:twisting-diff}
\w(d_+)=d_-,
\end{align}
see
\cite{FKW}
for the detail.

Assume $\kappa\ne 0$.
It was shown in \cite{FKW}
that
$H^{\bullet}_{\pm}(V)$,
$V\in \Obj \BGG_{\kappa}$,
is an object of $\Wcat$.
Thus,
we have a family of functors 
defined by
\begin{align}
 \begin{array}{cccc}
 \BGG_{\kappa} &\longrightarrow  &\Wcat& \\
V&\longmapsto & H^i_{\pm}(V)&(i\in \Z).
 \end{array}
\end{align}
\section{The 
Zhu algebra of $\W$-algebras}
\label{section:Zhu}
%
%
\subsection{The first isomorphism}
There is an obvious  homomorphism
$\U(\W_{\kappa}(\sg))
\rightarrow H^0(\U(C_{\kappa}(\g)))$,
which
induces
an  algebra homomorphism
$
\A(\W_{\kappa}(\sg))
\rightarrow H^0(\A(C_{\kappa}(\g)))
$.
On the other hand,
by \eqref{eq:iso:Zhu-tensor-1} and
\eqref{eq:realization_of_W_by_FB},
we have a commutative diagram
\begin{align}\label{eq:CD1}
\begin{CD}
\A(\W_{\kappa}(\sg)) @= \A(H^0( C_{\kappa}(\g)_0)))\\
@V{}VV @VV{}V\\
H^0(\A(C_{\kappa}(\g)))@= H^0(\A(C_{\kappa}(\g)_0)).
\end{CD}
\end{align}
But the vertical map on the right-hand-side
induces a map
$\A(\gr^G H^0( C_{\kappa}(\sg)_0))
\rightarrow \gr^G (H^0(\A(C_{\kappa}(\sg))))$,
which is an isomorphism 
by \eqref{eq-gr-zhu-1}
and 
\eqref{eq-gr-zhu2}.
Hence 
the map $\A(H^0( C_{\kappa}(\g)_0)))
\rightarrow H^0(\A(C_{\kappa}(\g)_0))$
is an isomorphism.
Therefore  we conclude as follows.
\begin{Pro}\label{Pro:zhu2}
For any complex number $\kappa$,
 \begin{align*}
  \text{$\A(\W_{\kappa}(\sg))
\cong H^0(\A(C_{\kappa}(\g)))
= H^0(\A(C_{\kappa}(\sg)_0))$
}
 \end{align*}as $\C$-algebras.
\end{Pro}

\subsection{The second isomorphism}
 Let
\begin{align*}
\bar{d}&=\bar{d}^{\st}+\bar{\chi}\in
U(\sg)\*\bar{\Cl} ,\\
\text{where }&\bar{d}^{\st}=\sum_{\alpha\in
\sproots}
J_{\alpha}(0){\psi}_{-\alpha}(0)
-\frac{1}{2}\sum_{{\alpha,\beta,\gamma\in
\sproots}{ }}
c_{\alpha,\beta}^{\gamma}
{\psi}_{-\alpha}(0){\psi}_{-\beta}(0)
\psi_{\gamma}(0),
\\
&\bar{\chi}=\sum_{\alpha\in \sPi}\psi_{-\alpha}(0).
\end{align*}
Then,
$(\bar{d}^{\st})^2=(\bar{\chi})^2=0$,
$\{\bar{d}^{\st},
\bar{\chi}\}=0$.
Thus,
$(\ad\bar{d})^2=0$
on $U(\sg)\*\bar{\Cl}$.
Here,
$U(\sg)\* \bar{\Cl}
$
is regarded as a superalgebra
and $\ad \bar{d}(u)$,
$u\in U(\sg)\* \bar{\Cl}$,
is the adjoint action on $u$ of the odd element $\bar{d}$.
Put
$\deg \psi_{\alpha}(0)=-1$,
$\deg \dual{\psi}_{\alpha}(0)=1$ ($\alpha\in 
\sproots$)
and $\deg u=0$ ($u\in U(\sg)$).
This makes 
$(U(\sg)\*\bar{\Cl},\ad \bar{d})$
a complex.
Note that the corresponding cohomology
$H^{\bullet}(U(\sg)\*\bar{\Cl},\ad \bar{d})$
is naturally a graded $\C$-algebra.
The following proposition is straightforward
from Proposition
\ref{Pro:Zhu_of_complex}. 
\begin{Pro} \label{Pro:iso-pre-Kostant}
We have the natural isomorphism of graded $\C$-algebras 
\begin{align*}
H^{\bullet}(\A(C_{\kappa}(\g))
)\cong H^{\bullet}(U(\sg)\*\bar{\Cl},\ad \bar{d}).
\end{align*}
\end{Pro}
By \eqref{eq:iso:Zhu-tensor-1}
and Proposition \ref{Pro:iso-pre-Kostant},
it follows that
\begin{align}
 H^i(U(\sg)\*\bar{\Cl},\ad \bar{d})=\{0\}
\quad (i\ne 0).
\end{align}
\begin{Pro}\label{Pro:zhu1}
The map defined by
\begin{align*}
\begin{array}
{ccc}
\Center(\sg)&\rightarrow &H^0(U(\sg)\*\bar{\Cl},\ad \bar{d})\\
z&\mapsto &z \*1
\end{array}
\end{align*}
is an isomorphism of
 $\C$-algebras.
\end{Pro}
\begin{proof}
{\textbf{Step 1)}}\
In what follows we shall identify
$\Lam (\sn_-)$ with $\Lam (\dual{\sn_+})$.
Thus,
$\bar{\Cl}=\Lam (\dual{\sn_+})\*
\Lam (\sn_+)$
as a $\C$-vector space.
Set
$\psi_{\alpha}=\psi_{\alpha}(0)\in
\Lam (\sn_+)$
and $\dual{\psi}_{\alpha}=\psi_{-\alpha}(0)\in
\Lam (\dual{\sn_+})$.
Put
\begin{align*}
\bar{C}=
\sum_{n\in \Z}\bar{C}^n=U(\sg)\*\bar{\Cl}.
\end{align*}
Here,
\begin{align}
\bar{C}^n=\sum_{i-j
=n} U(\sg)\* \Lam^i (\dual{\sn_+})
\*
\Lam^j({\sn_+}).
\end{align}
We have
\begin{align}\label{eq:dec_bar_d}
\ad \bar{d}=\bar{d}_++\bar{d}_-
\text{ }
\end{align}
on $\bar{C}$,
where
$\bar{d}_{\pm}$
is
defined by 
\begin{equation}\label{eq:bar-d-+}
\begin{aligned}
\bar{d}_+(u\*\omega_1\*
\omega_2)
=\sum_{\alpha
\in \sproots}
\left\{(
(\ad J_{\alpha}
(u))
\* \dual{\psi}_{\alpha}\omega_1\*
\omega_2
+u\*\dual{\psi}_{\alpha}\omega_1\*
\ad J_{\alpha}(\omega_2)
\right\}\\
-\frac{1}{2}\sum_{\alpha,
\beta,\gamma\in \sproots}
c_{\alpha,\beta}^{\gamma} u\*
\dual{\psi}_{\alpha}\dual{\psi}_{\beta}
(\ad \psi_{\gamma}(\omega_1))
\* \omega_2,
\end{aligned}
\end{equation}
\begin{equation}\label{eq:bar-d--}
\begin{aligned}
(-1)^{i}\bar{d}_-(u\*\omega_1\*
\omega_2)
=\sum_{\alpha\in \sproots}
u(J_{\alpha}+\bar{\chi}(J_{\alpha}))\*
\omega_1\*
\ad \dual{\psi}_{\alpha}
(\omega_2)
\\
-\frac{1}{2}\sum_{\alpha,
\beta,\gamma\in \sproots}
c_{\alpha,\beta}^{\gamma} u\*
\omega_1\*\psi_{\gamma}\ad \dual{\psi}_{\beta}
(\ad 
\dual{\psi}_{\alpha}
(\omega_2))
\end{aligned}
\end{equation}
for $u\in U(\sg)$,
$\omega_1\in \Lam^i(\sn_+^*)$,
$\omega_2\in \Lam(\sn_+)$.
Note that
\begin{equation}\label{eq:bar-d+-d-}
\begin{aligned}
\bar{d}_+(U(\sg)\* \Lam^i(\dual{\sn_+})\*
\Lam^{j}(\sn_+))
\subset
U(\sg)\* \Lam^{i+1}(\dual{\sn_+})\*
\Lam^{j}(\sn_+),
\\
\bar{d}_-(U(\sg)\* \Lam^i(\dual{\sn_+})\*
\Lam^{j}(\sn_+))
\subset
U(\sg)\* \Lam^{i}(\dual{\sn_+})\*
\Lam^{j-1}(\sn_+).
\end{aligned}
\end{equation}
Hence,
it follows that
\begin{align}\label{ew:commutativity-dif-bar}
\bar{d}_+^2=\bar{d}_-^2=\{\bar{d}_+,\bar{d}_-\}=0.
\end{align}

Set
\begin{align}\label{eq:filtration-Kostant}
F^p \bar{C}=\sum_{i\geq p}U(\sg)\* \Lam^i(\dual{\sn_+})\*
\Lam^{\bullet}(\sn_+)\subset \bar{C},
\end{align}
Then,
\begin{align}
&\bar{C}=F^0\bar{C}
\supset F^1\bar{C}\supset \dots
\supset F^{\dim \sn_+ +1}\bar{C}=\{0\},
\\&\bigcap F^p \bar{C}=\{0\},\\
&
d_+ F^p \bar{C}
\subset F^{p+1} \bar{C},\quad
d_- F^p \bar{C}\subset F^{p} \bar{C}.
\end{align}
%
%
%
%
%
%
We  consider the  spectral sequence
$E$
such that
\begin{align*}
 E_0^{p,q}=F^p \bar{C}^{p+q}/F^{p+1}\bar{C}^{p+q}.
\end{align*}

\medskip

\noindent {\textbf{Step 2)}}\
We have
$E_1
=H^{\bullet}(\bar{C},\bar{d}_-)$. 
Let
$\C_{\bar{\chi}}$ be the one-dimensional
representation of $U(\sn_+)
$ defined by the character
$\bar{\chi}$.
Then,
 by
\eqref{eq:bar-d--},
we see that  
the complex
$(\bar{C},\bar{d}_-)$
is nothing but the Chevalley complex
for calculating the $\sn_+$-homology
$H_{\bullet}(\sn_+,\left(U(\sg)\* \C_{\bar{\chi}}\right)\*
\Lam (\dual{\sn}_+))$
(with the opposite grading).
Here,
$\left(U(\sg)\* \C_{\bar{\chi}}\right)\*
\Lam (\dual{\sn}_+)$
is regarded as a right $U(\sn_+)$-module
on which
  $U(\sn_+)$
acts on the first factor
$U(\sg)\* \C_{\bar{\chi}}$.
But obviously,
$U(\sg)$
is free over $\sn_+$,
thus, so is $U(\sg)\* \C_{\bar{\chi}}$.
Therefore,
\begin{equation}
 \begin{aligned}
E_1^{\bullet,q}&=\begin{cases}
\left[\left(U(\sg)\* \C_{\bar{\chi}}\right)
/\left(U(\sg)\* \C_{\bar{\chi}}\right)\sn_+
\right]
\*\Lam (\dual{\sn}_+)
&(q=0)\\
0&(q\ne 0)
\end{cases}\\
&=\begin{cases}
\left(U(\sg)\*_{U(\sn_+)}\C_{-\bar{\chi}}
\right)
\*\Lam (\dual{\sn}_+)
&(q=0)\\
0&(q\ne 0).
\end{cases}
\end{aligned}\label{eq:E_1-Kostant}
\end{equation}
Note that the space  $
U(\sg)\*_{U(\sn_+)}\C_{-\bar{\chi}}
$
 has   a natural left
$\sn_+$-module structure defined 
by 
$x\cdot u\* 1=(x u)\*1$,
$x\in \sn_+, u\in U(\sg)$.

\medskip

\noindent {\textbf{Step 3)}}\
Next we calculate
$E_2=H^{\bullet}(H^{\bullet}(\bar{C},\bar{d}_-),\bar{d}_+)$.
The formula
 \eqref{eq:bar-d-+}
and \eqref{eq:E_1-Kostant}
 show 
that
\begin{align*}
E_2^{p,0}=H^{p}(\sn_+,
\left(U(\sg)\*_{U(\sn_+)}\C_{-\bar{\chi}}\right)
\* \C_{\bar{\chi}}).
\end{align*}
Here,
the action of $\sn_+$
on $\left(U(\sg)\*_{U(\sn_+)}\C_{-\bar{\chi}}\right)
\* \C_{\bar{\chi}}$ is the tensor product
action.
But
it is well-known  since Kostant 
\cite{Kostant-Whittaker} that
$H^{i}(\sn_+,
\left(U(\sg)\*_{U(\sn_+)}\C_{-\bar{\chi}}\right)
\* \C_{\bar{\chi}})=0$
($i\ne 0$)
and
we have an isomorphism
\begin{align*}
\begin{array}
{ccc}
\Center(\sg)&\cong &E_2^{0,0}=H^{0}(\sn_+,
\left(U(\sg)\*_{U(\sn_+)}\C_{-\bar{\chi}}\right)
\* \C_{\bar{\chi}}),\\
z&\mapsto &(z\* 1)\* 1.
\end{array}
\end{align*}
Hence the spectral sequence collapse at $E_2=E_{\infty}$
and Proposition is proved.
\end{proof}
By Proposition \ref{Pro:zhu2}
and Proposition \ref{Pro:zhu1},
we conclude as the following.
\begin{Th}\label{Th:Zhu-Center}
For any complex number $\kappa$,
the Zhu algebra
$\A(\W_{\kappa}(\sg))$
of $\W_{\kappa}(\sg)$
is naturally isomorphic to
$\Center(\sg)$.
\end{Th}
\begin{Rem}\label{Rem:image.of.L(0)}
Let $\kappa\ne 0$.
Then,
one calculates that,
 under the isomorphism
$\A(\W_{\kappa}(\sg))\cong \Center(\sg)$,
the image of $L(0)$ in $\A(\W_{\kappa}(\sg))$
is mapped to the element
\begin{align*}
 \frac{1}{2\kappa}\Omega
-\frac{1}{2}
\left(
{\kappa}|\srho\che|^2-2\bra\srho,\srho\che\ket
\right)\id
\end{align*} 
of $\Center(\sg)$.
Here,
$\Omega$ is the Casimir element of $U(\sg)$.
\end{Rem}

\section{Simple Modules over $\W$-algebras}
\subsection{Verma modules of 
$\W$-algebra}\label{subsection:Verma}
In what follows we identify
$\A(\W_{\kappa}(\sg))$
with $\Center(\sg)$
by 
 Theorem \ref{Th:Zhu-Center}.
Let
$\gamma: \Center(\sg)\cong S(\sh)^{\sW}$
be the Harish-Chandra isomorphism.
For $\bar{\lam}\in \dual{\sh}$,
let
\begin{align*}
\gamma_{\bar{\lam}}=\text{(evaluation at $\bar{\lam}+\srho$)}
\circ{\gamma}:
\Center(\sg)\rightarrow \C.
\end{align*}
Thus,
$z\in \Center(\sg)$
acts as $\gamma_{\bar{\lam}}(z)\id$
on the Verma module of $\sg$
of highest weight $\bar{\lam}$.
Let $\C_{\gamma_{\bar{\lam}}}$
 be
one-dimensional
representation of $\Center(\sg)=\A(\W_{\kappa}(\sg))$
defined by $\gamma_{\bar{\lam}}$.
We also regard $\C_{\gamma_{\bar{\lam}}}$  
 as a 
$\U(\W_{\kappa}(\sg))_{\geq 0}$-module
on which
$\U(\W_{\kappa}(\sg))_n$,
$n>0$,
acts trivially.
Define
\begin{align}
\M(\gamma_{\bar{\lam}})=\U(\W_{\kappa}(\sg))\*_{
\U(\W_{\kappa}(\sg))_{\geq 0}}\C_{\gamma_{\bar{\lam}}}.
\end{align}
The module $\M(\gamma_{\bar{\lam}})$
is called the {\em{Verma module of $\W_{\kappa}(\sg)$
of highest weight $\gamma_{\bar{\lam}}$}}.
Let $|\gamma_{\bar{\lam}}\ket $
denote the vector  $1\* 1\in \M(\gamma_{\bar{\lam}})$.
By \eqref{eq:cc} and Remark \ref{Rem:image.of.L(0)},
we have:
\begin{align}
 L(0)|\gamma_{\bar{\lam}}\ket 
=\Delta_{\bar{\lam}}|\gamma_{\bar{\lam}}\ket \quad \text{(if $\kappa\ne0$)},
\end{align}
where
\begin{align}
  \Delta_{\bar{\lam}}
=\frac{|\bar{\lam}+\srho|^2}{2\kappa}
-\frac{\rank \sg}{24}+\frac{c(\kappa)}{24}.
 \end{align}

 \smallskip

Let 
\begin{align*}
 \text{$G^p \M(\gamma_{\bar{\lam}})=\sum_{p_1+\dots +p_r\geq p}
G^{p_1}\W_{\kappa}(\sg)\dots G^{p_r}\W_{\kappa}(\sg)|\gamma_{\bar{\lam}}
\ket
\subset \M(\gamma_{\bar{\lam}})$
\quad ($p\leq 1$).
}
\end{align*}Then,
\begin{align}
 & \dots \supset G^p\M(\gamma_{\bar{\lam}})
\supset G^{p+1}\M(\gamma_{\bar{\lam}})\supset 
\dots \supset G^0\M(\gamma_{\bar{\lam}})
\supset G^1\M(\gamma_{\bar{\lam}})= \{0\},
\\
&\M(\gamma_{\bar{\lam}})= \bigcup_p G^p \M(\gamma(\bar{\lam}))
\\& 
G^p \W_{\kappa}(\sg) \cdot G^q \M(\gamma_{\bar{\lam}})
\subset G^{p+q}\M(\gamma_{\bar{\lam}}).\label{eq:graded-action}
\end{align}
Here,
in \eqref{eq:graded-action},
$G^p \W_{\kappa}(\sg)
 \cdot G^q \M(\gamma_{\bar{\lam}})$
denotes the span of the vectors
$Y_n(v)m$
($n\in \Z$,
$v\in G^p \W_{\kappa}(\sg)$
 $m\in G^q \M(\gamma_{\bar{\lam}})$).
Let $\gr^G \M(\gamma_{\bar{\lam}})$
be the corresponding graded 
vector space.
By \eqref{eq:graded-action},
$\gr^G \M(\gamma_{\bar{\lam}})$ is a module over 
the commutative vertex algebra $\gr^G \W_{\kappa}(\sg)
$.
Let $\bar{W}_i(n)$
denote the image of $W_i(n)$
in $\gr^G \W_{\kappa}(\sg)$.
The following proposition is easy to see.
\begin{Pro}
For all $\bar{\lam}\in \dual{\sh}$,
we have the isomorphism
\begin{align*}
\begin{array}{ccc}
\C [\bar{W}_1(-n_1),\dots , \bar{W}_{l}(-n_l)]_{n_1,\dots ,n_l\geq 1}
&\cong &
\gr^G \M(\gamma_{\bar{\lam}})\\
f&\mapsto &f \overline{|\gamma_{\bar{\lam}}\ket }
\end{array}
\end{align*}
Here, $\overline{|\gamma_{\bar{\lam}}\ket }$
is the image of ${|\gamma_{\bar{\lam}}\ket }$
in $\gr^G \M(\gamma_{\bar{\lam}})$.
\end{Pro}

\begin{Co}\label{Co:basis-of-Verma}
For  $\bar{\lam}\in \dual{\sh}$,
the set
\begin{align*}
 \left\{
W_{i_1}(-n_1)\dots
W_{i_m}(-n_m)|\gamma_{\bar{\lam}}\ket ;
\begin{array}
{l}
1\leq i_1\leq \dots \leq i_m\leq l,
n_p\geq 1,
\\n_p\geq n_{p+1}\text{ if
$i_p=i_{p+1}$ }
\end{array}\right\}.
\end{align*}
forms a basis of $\M(\gamma_{\bar{\lam}})$.
In particular,
for $\kappa \in \C^*$,
$\M(\gamma_{\bar{\lam}})\in \Wcat$
and 
\begin{align*}
\ch \M(\gamma_{\bar{\lam}})=\frac{q^{
\frac{|\bar{\lam}+\srho|^2}{2\kappa}
}}{
\eta(\tau)^{\rank\sg}}.
\end{align*}
where $
\eta(\tau)=q^{\frac{1}{24}}\prod_{i\geq 1}(1-q^i)$,
$q=e^{2\pi\sqrt{-1}\tau}$.
\end{Co}

Assume $
\kappa \ne 0$.
By Corollary \ref{Co:basis-of-Verma},
we have 
\begin{align}\label{eq:hw-of-Verma}
 \M(\gamma_{\bar{\lam}})
=\bigoplus_{\Delta\in \Delta_{\bar{\lam}}+
\Z_{\geq 0}}\M(\gamma_{\bar{\lam}})_{\Delta},
\quad \M(\gamma_{\bar{\lam}})_{\Delta_{\bar{\lam}}}=\C
 |\gamma_{\bar{\lam}}
\ket.
\end{align}
Therefore,
 $N_{\Delta_{\bar{\lam}}}=\{0\}$
for any proper submodule $N$ of $\M(\gamma_{\bar{\lam}}$).
Hence,
the Verma module
$\M(\gamma_{\bar{\lam}})$
has a unique simple quotient
which we shall denote
by $\why(\gamma_{\bar{\lam}})$.
The module
$\why(\gamma_{\bar{\lam}})$
is 
called the {\em irreducible $\W_{\kappa}(\sg)$-module of highest
weight $\gamma_{\bar{\lam}}$}.
The following theorem is clear
by \cite[Theorem 1.4.2]{FZ}.
\begin{Th}
For any $\kappa\in \C^*$,
the set
\begin{align*}
\{\why(\gamma_{\bar{\lam}});\bar{\lam}\in
\dual{\sh}/\sim\}
\end{align*}
is a complete set of isomorphism classes of simple
modules in $\Wcat$.
Here,
$\sim$ is an equivalence relation
defined by $\bar{\lam}\sim \bar{\mu}\iff
\bar{\mu}+\srho\in \sW(\bar{\lam}+\srho)$
$(\bar{\lam},\ \bar{\mu}\in \dual{\sh})$.
\end{Th}

\subsection{The duality functor $D$}
Assume that $\kappa \ne 0$
so that $L(z)$ is well-defined.

For  $V\in \Wcat$
(or for simply a $L(0)$-module $V$),
let
\begin{align}\label{eq:def-of-D(V)}
 \text{$D(V)=\bigoplus_{h\in \C}\Hom_{\C}(V_{h}^{\gen},\C)$.}
\end{align}
Let
\begin{align}
\theta (Y_n(v))=(-1)^{\Delta}\sum_{j\geq 0}
\frac{1}{j!}Y_{-n}(L(1)^jv)\quad
\text{for $v\in \W_{\kappa}(\sg)_{\Delta}$ and $n\in
\Z$.}
\end{align}
Then,
as is known,
the space
$D(V)$
has a $\W_{\kappa}(\sg)$-module structure
defined
by
$\bra Y_n(v)f,v_1\ket=\bra f, \theta(Y_n(v))v_1\ket$
 ($v\in \W_{\kappa}(\sg)$,
$f\in D(V)$,
$v_1\in V$).
The correspondence $V\fmap D(V)$
defines a duality functor in $\Wcat$.

It is known that
the map $\theta$ induces
an
 anti-automorphism
of 
$\A(\W_{\kappa}(\sg))$,
which is also  denoted by $\theta$.

\begin{Lem}\label{Lem:anti-auto-on-Zhu}
Let $\bar{\theta}$
be the anti-automorphism of $U(\sg)$
defined by
\begin{align*}
&\bar{\theta}(J_{\pm \alpha})=
-(-1)^{\height \alpha}J_{\pm \alpha}~
(\alpha\in \sproots),\quad
\bar{\theta}(h)=-h~(h\in \sh).
\end{align*}
Then,
$\theta(z)=\bar{\theta}(z)$
for $z\in \Center(\sg)$
under the identification 
$\A(\W_{\kappa}(\sg))=
\Center(\sg)$.
\end{Lem}
\begin{proof}
Straightforward  from the definition.
\end{proof}
\begin{Lem}\label{Lem:dual-center}
For $\bar{\lam}\in \dual{\sh}$,
$\gamma_{\bar{\lam}}\circ
\bar{\theta}_{\mid \Center(\sg)}=\gamma_{-w_0(\bar{\lam})}$.
\end{Lem}
\begin{proof}
Since
$-w_0(\lam)+\srho=-w_0(\lam+\srho)$,
it is suffucient
to show that
\begin{align}\label{eq:inv_and_harish-chandra}
(\gamma\circ \bar{\theta}(z))(\lam)=\gamma(z)(-\lam)
\quad\text{(for $z\in \Center(\sg)$)}
\end{align}
for all $\lam\in \dual{\sh}$.
But certainly  \eqref{eq:inv_and_harish-chandra}
holds for all
$\lam\in \sP_++\srho$,
where
$\sP_+=\{ \lam\in \sP;
\bra \lam,\alpha\che\ket\geq 0\text{ for all }
\alpha\in \sproots\}\}$.
Therefore,
\eqref{eq:inv_and_harish-chandra}
holds for all $\lam\in \dual{\sh}$.
%
\end{proof}
Lemma \ref{Lem:anti-auto-on-Zhu}
and
Lemma \ref{Lem:dual-center}
imply the following result.
\begin{Th}\label{Th:Duality}
For all $\bar{\lam}\in \dual{\sh}$,
we have
$D(\why(\gamma_{\bar{\lam}}))\cong \why(\gamma_{-w_0(\bar{\lam})})$.
\end{Th}

The functor $D$
is clearly exact.
Thus,
by Theorem \ref{Th:Duality},
the exact sequence $\M(\gamma_{-w_0(\bar{\lam})})
\rightarrow
\why(\gamma_{-w_0(\bar{\lam})})
\rightarrow 0$
gives rise to the exact sequence 
\begin{align}
 0
\rightarrow \why(\gamma_{\bar{\lam}})
\rightarrow D(\M(\gamma_{-w_0(\bar{\lam})})).
\end{align}
Thus,
$D(\M(\gamma_{-w_0(\bar{\lam})}))$ has $\why(\gamma_{\bar{\lam}})$
as its unique simple submodule.
The following lemma will be used in next sections.
\begin{Lem}\label{Lem:criterion-simple}
 Let $V$ be an object of $\Wcat$.
Suppose there exist exact sequences
\begin{align}
&\M(\gamma_{\bar{\lam}})
\rightarrow
V
\rightarrow 0,
\label{eq:surj-Lem} \\
 &0
\rightarrow V
\rightarrow D(\M(\gamma_{-w_0(\bar{\lam})})).
\label{eq:inj-Lem} 
\end{align}
Then,
$V$ is either $\{0\}$ or isomorphic to $\why(\gamma_{\bar{\lam}})$.
\end{Lem}
\begin{proof}
 Let
$N$ be a proper submodule of $V$.
Then,
\eqref{eq:surj-Lem} 
implies that
$N_{\Delta_{\bar{\lam}} }=\{0\}$.
But  then,
 by \eqref{eq:inj-Lem},
it follows that
$N=\{0\}$.
Lemma is proved.
\end{proof}
Now,
for $\bar{\lam}\in \dual{\sh}$,
set
\begin{align*}
 \text{$G^p D(\M(\gamma_{\bar{\lam}}))
=D(\M(\gamma_{\bar{\lam}})/G^{-p}\M(\gamma_{\bar{\lam}}) )$
\quad (as a vector space).
}
\end{align*}Then,
\begin{align*}
 &D(\M(\gamma_{\bar{\lam}}))
=
G^{-1} D(\M(\gamma_{\bar{\lam}}))\supset 
G^0  D(\M(\gamma_{\bar{\lam}}))\supset \dots
\supset 
 G^p D(\M(\gamma_{\bar{\lam}}))\supset \dots, \\
& \bigcap_p G^p D(\M(\gamma_{\bar{\lam}}))=\{0\},\\
& G^p \W_{\kappa}(\sg)\cdot G^q D(\M(\gamma_{\bar{\lam}}))
\subset G^{p+q}D(\M(\gamma_{\bar{\lam}})).
\end{align*}
Let $\gr^G D(\M(\gamma_{\bar{\lam}}))$
be the corresponding graded vector space.  
Then,
\begin{align*}
\gr^G D(\M(\gamma_{\bar{\lam}}))=
\left( \gr^G \M(\gamma_{\bar{\lam}})\right)^*,
\end{align*} 
where ${}^*$ denotes the graded dual in the obvious sense.
The space $\gr^G D(\M(\gamma_{\bar{\lam}}))
$
 is a module over $\gr^G \W_{\kappa}(\sg)$,
on which $\overline{W}_i(n)$ ($i \in \sI$, $n\in \Z$)
acts as 
\begin{align}
(\overline{W}_i(n)f)(v)=-(-1)^{d_i}f(\overline{W}_i(-n)v)\quad
(f\in \gr^G D(\M(\gamma_{\bar{\lam}})),
v \in \gr^G \M(\gamma_{\bar{\lam}})).
\end{align}
In particular,
it follows that
\begin{align}\label{eq:singular-vector-of-dual-of-Verma}
\{v\in \gr^G D(\M(\gamma_{\bar{\lam}}));
\bar{W}_i(n)v=0\quad (i\in \sI, n>0)\}
=
\C \overline{|\gamma_{\bar{\lam}}}\dual{\ket}.
\end{align}
Here,
$\overline{|\gamma_{\bar{\lam}}}\dual{\ket}$
is the image of the vector $|\gamma_{\bar{\lam}}\dual{\ket}$
dual to $|\gamma_{\bar{\lam}}\ket$.

\section{Quantized reductions
and the representation theory of $\W$-algebras}
In this section we assume that
$\kappa\in \C^*$ unless otherwise stated.
\subsection{The images of Verma modules}
For an object $V$ of $\Wcat$,
let 
\begin{align}
\text{$V_{\tp}=\sum\limits_{h\in \C\atop
V_{h-n}^{\gen}=\{0\}\text{ for all $n>0$}} V_h^{\gen}\subset V$.}
\end{align}
Then,
$\U(\W_{\kappa}(\sg))_n \cdot V_{\tp}=\{0\}$
for $n>0$.
Thus,
$V_{\tp}$
is naturally an $\A(\W_{\kappa}(\sg))$-module.

By \cite[Theorem 5.7, Remark 5.8]{A},
we have:
\begin{align}
&H_{\pm}^i(M(\lam))=\{0\} \quad (i\ne 0),\\
&
 \ch H_-^0(M(\lam))=\ch \M(\gamma_{\bar{\lam}}),
\label{eq:ch-Verma-}\\
&
\ch H_+^0(M(\lam))=\ch \M(\gamma_{\overline{t_{-\srho\che}\circ \lam}}),
\label{eq:ch-Verma+}
\end{align}
for all $\lam\in \dual{\h}$.
 By \eqref{eq:ch-Verma+},
one sees that
\begin{align}
& H_-^0(M(\lam))_{\tp}
=H_-^0(M(\lam))_{\Delta_{\bar{\lam}}}
=\C |\lam\ket,\\&
 H_+^0(M(\lam))_{\tp}
=H_+^0(M(\lam))_{\Delta_{\overline{t_{-\srho\che}\circ \lam}}}
=\C |\lam\ket.
\end{align}
Here,
$|\lam\ket =v_{\lam}\* \1$
and $v_{\lam}$
is the highest weight vector of $M(\lam)$.
\begin{Lem}\label{Lem:hw-Verma--}
$ $

\begin{enumerate}
 \item For any $\lam\in \dual{\h}$,
$ H_-^0(M(\lam))_{\tp}
\cong \C_{\gamma_{\bar{\lam}}}$
as $\A(\W_{\kappa}(\sg))$-modules.
 \item For any $\lam\in \dual{\h}$,
$ H_+^0(M(\lam))_{\tp}
\cong \C_{\gamma_{\overline{t_{-\srho\che}\circ \lam}}}$
as $\A(\W_{\kappa}(\sg))$-modules.
\end{enumerate}\end{Lem}
\begin{proof}
(1) It is not difficult to see that 
the center
$\Center(\sg)$
naturally acts on the 
space $H_-^0(M(\lam))_{\tp}
=\C |\lam\ket $
and that
its action coincides with the action
of $\A(\W_{\kappa}(\sg))$ 
under the identification 
$\A(\W_{\kappa}(\sg))=\Center(\sg)$.
(2)
Let
$|\lam\ket_{\pm}$
temporary denote
the vector $v_{\lam}\* \1\in M(\lam)\* \F(\Ln{\pm})$.
Since $\A(\W_{\kappa}(\sg))=H^0(\A(C_{\kappa}(\sg)_0))
$
by Proposition \ref{Pro:zhu2},
it is sufficient to show that
$\C |\lam \ket_+\cong \C |\w \circ \lam\ket_-$
as $\A(C_{\kappa}(\sg)_0)$-modules.
Here,
$\w$ is defined in Section \ref{section:functors}.
But $\A(C_{\kappa}(\sg)_0)$
is generated by
$\widehat{J}_i(0)$ ($i\in \sI$),
$\widehat{J}_{-\alpha}(\height \alpha)$,
$\psi_{-\alpha}(\height \alpha)$
$(\alpha \in \sproots)$,
and
both spaces
$\C |\lam\ket_{+} $
and $\C |\lam\ket_{-} $
are   $\A(C_{\kappa}(\sg)_0)$-modules on which 
$\widehat{J}_{-\alpha}(\height \alpha)$,
$\psi_{-\alpha}(\height \alpha)$
$(\alpha \in \sproots)$
act trivially.
Hence,
$\C |\lam\ket_+\cong \C |\w \circ \lam \ket_-$
as $\A(\W_{\kappa}(\sg)_0)$-modules by \cite[Proposition 4.2]{A}.
\end{proof}

Now we have the following result.
\begin{Th}\label{Th:image_of_Verma_Module}
For any $\lam\in \dual{\h}$,
\begin{align*}
H_-^0(M(\lam))\cong \M(\gamma_{\bar{\lam}}).
\end{align*}
\end{Th}
\begin{proof}
By Lemma \ref{Lem:hw-Verma--} (1),
we have a homomorphism
$\phi:
\M(\gamma_{\bar{\lam}})\rightarrow \F(M(\lam))$
of $\W_{\kappa}(\sg)$-modules
define by $\phi(|\gamma_{\bar{\lam}}\ket)=
|\lam \ket$.
Hence,
 it is sufficient to show that
$\phi$ is a bijection.
This follows from the following claim:

\medskip

\noindent {\em{Claim}} \
There exists a decreasing filtration 
$\{G^p H_-^0(M(\lam))\}_{p\leq 1}$
of $H_-^0(M(\lam))$
such that
\begin{align*}
& \text{$\dots \supset G^p H_-^0(M(\lam))\supset
G^{p+1}H_-^0(M(\lam))\supset \dots
\supset
G^1 H_-^0(M(\lam))=\{0\}$}, 
\\ & H_-^0(M(\lam))=\bigcup_{p\leq 0}G^pH_-^0(M(\lam)),\\
&\phi(G^p(\M(\gamma_{\bar{\lam}})))\subset G^p H_-^0(M(\lam))
\quad (\forall p),
\end{align*}
and 
the induced map
$\bar{\phi}: \gr^G \M(\gamma_{\bar{\lam}})\rightarrow \gr^G
 H_-^0(M(\lam))$
is an isomorphism of $\gr^G \W_{\kappa}(\sg)$-modules.

\medskip 

The above claim can be seen using the argument of \cite{FB}
as follows:
Let $C(\lam)_0$
be  the subspace of
$C(\Ln{-},M(\lam))$ spanned by the elements
\begin{align*}
\widehat{J}_{a_1}(-n_1)\dots
\widehat{J}_{a_p}(-n_p)
{\psi}_{\alpha_{j_1}}
(-m_1)\dots {\psi}_{\alpha_{j_q}}(-m_q)
|\lam \ket\quad
(a_i\in \sproots\sqcup \sI,
\alpha_{j_i}\in \sproots).
\end{align*}
Then,
$d_-C(\lam)_0\subset C(\lam)_0$
and
\begin{align}\label{eq:coho-Verma}
 \text{$H_-^{\bullet}(M(\lam))=
H^{\bullet}(C(\lam)_0)
$,
}\quad H^{i}(C(\lam)_0)=0\ (i\ne 0),
\end{align}see \cite[section 5]{A}.
Define a decreasing filtration $\{G^p C^n(\lam)_0\}$
on $C(\lam)_0$
by setting 
\begin{align*}
G^p C^n(\lam)_0
=
\bigoplus_{\mu
\atop \bra\mu-\lam,\srho\che\ket
\leq n-p}C^{n}(\lam)_0^{\mu}
\end{align*}
Then,
one can check that 
the corresponding spectral sequence converges,
and
it induces
a filtration on $H_-^0(M(\lam))$
with desired properties.
\end{proof}
\begin{Rem}
\begin{enumerate}
 \item  It is not likely that
$H_+^0(M(\lam))\cong 
\M(\gamma_{\overline{t_{-\srho\che}\circ {{\lam}}}})
$ in general.
\item By replacing $L(0)$
with the formal degree operator defined
in \cite{A},
one sees that Theorem \ref{Th:image_of_Verma_Module}
holds also for $\kappa=0$.
\end{enumerate}
\end{Rem}
\subsection{The images of the duals of Verma modules}
For $V\in \Obj\Wcat$,
let
\begin{align*}
\HW(V)=\{v\in V;
W_i(n)v=0
\text{ for all $i=1,\dots,l$ and $n>0$}\}.
\end{align*}
It is naturally
an
$\A(\W_{\kappa}(\sg))$-module.
Clearly,
$V_{\tp}\subset \HW(V)$.
\begin{Lem}\label{Lem:dual-of-Verma}
Let $\bar{\lam}\in \dual{\sh}$.
Let $M$
be an object in
$\Wcat$
such that
 $\ch M=\ch \M(\gamma_{\bar{\lam}})$
and $\HW (M)\cong \C_{\gamma_{\bar{\lam}}}$
as $\A(\W_{\kappa}(\sg))$-modules.
Then,
$M\cong D(\M(\gamma_{-w_0(\bar{\lam})}))$.
\end{Lem}
\begin{proof}
By the assumption,
$M=\bigoplus\limits_{\Delta
\in \Delta_{\bar{\lam}}+\Z_{\geq 0}}M_{\Delta}$
and $M_{\Delta_{\bar{\lam}}}\cong
\C_{\gamma_{\bar{\lam}} }$.
Thus,
$D(M)=\bigoplus\limits_{\Delta
\in \Delta_{\bar{\lam}}+\Z_{\geq 0}}D(M)_{\Delta}$
and $D(M)_{\Delta_{\bar{\lam}}}\cong
\C_{\gamma_{-w_0(\bar{\lam})}}$
by Lemma \ref{Lem:dual-center}.
Hence,
there exists a homomorphism
$\phi:
\M(\gamma_{-w_0(\bar{\lam})})
\rightarrow
D(M) $
in $\Wcat$.
Consider the exact sequence
$\M(\gamma_{-w_0(\bar{\lam})})
\rightarrow D(M)\rightarrow \coker \phi
 \rightarrow 0$.
It induces an exact sequence
\begin{align}
 \text{$0\rightarrow D(\coker \phi)\rightarrow M\rightarrow 
D(\M(\gamma_{-w_0(\bar{\lam})}))$.
}
\end{align}In particular,
$D(\coker \phi)\subset M$.
Thus,
$\HW(D(\coker \phi))
\subset \HW(M)=M_{\Delta_{\bar{\lam}}}$.
 This implies
$\HW(D(\coker \phi))=\{0\}$,
for
$D(\coker \phi)_{\Delta_{\bar{\lam}}}=0$.
But, then, $D(\coker \phi)=\{0\}$.
Therefore,
$M$ is a submodule of $D(\M(\gamma_{-w_0(\bar{\lam})}))$.
Hence,
the assumption 
$\ch M=\ch D(\M(\gamma_{-w_0(\bar{\lam})}))
$
shows that
 $M=D(\M(\gamma_{-w_0(\bar{\lam})}))$.
\end{proof}
\begin{Th}\label{Th:dual}
Let $\lam\in \dual{\h}$
be non-critical.

\begin{enumerate}
\item 
Suppose  that
$\bra \lam+\rho,\alpha\che\ket\not\in
\Z_{\geq 1}$
for all  $\alpha\in \sproots$.
Then,
$H_-^0(M(\lam)^*)\cong D(\M(\gamma_{-w_0(\bar{\lam})}))$.
\item 
Suppose  that
$\bra \lam+\rho,\alpha\che\ket\not\in \Z_{\geq
1}$ for
all $\alpha\in \prroots\cap t_{\srho\che}(\nrroots)$.
Then,
$H_+^0(M(\lam)^*)\cong D(\M(\gamma_{-w_0(\overline{
t_{-\srho\che}\circ \lam})}))$.
\end{enumerate}
\end{Th}
\begin{proof}
(1) By 
\cite[Theorem 6,8 (2)]{A},
we have
\begin{align}
 & H^i_-(M(\lam)^*)=\{0\}\quad (i\ne 0),\\
& H^0_-(M(\lam)^*)\cong D(\M(\gamma_{\bar{\lam}}))
\text{ as $L(0)$-modules.}
\label{eq:iso-dual-Verma-fact}
\end{align}
under the assumption of Theorem (1).
In particular,
\begin{align}\label{eq:ch-equal}
\text{$\ch H_-^0(M(\lam)^*)=
\ch \M(\gamma_{\bar{\lam}})$.
}\end{align}
This further implies that
 \begin{align}\label{eq:iso-as-Zhu}
\text{$H_{-}^0(M(\lam)^*)_{\Delta_{\bar{\lam}}}=\C |\lam\dual{\ket}
\cong 
\C_{\gamma_{\bar{\lam}}}$
as  $\A(\W_{\kappa}(\sg))$-modules.
}
 \end{align}
Here,
$|\lam\dual{\ket}
=\dual{v}_{\lam}
\* \1$
and $\dual{v}_{\lam}\in M(\lam)^*$
is the vector dual to the highest weight vector $v_{\lam}\in M(\lam)$.
Thus,
by Lemma \ref{Lem:dual-of-Verma},
it is sufficient to show 
that
$\HW( H_{-}^0(M(\lam)^*))=
\C |\lam\dual{\ket}$.

Now
the isomorphism 
in
\eqref{eq:iso-dual-Verma-fact}
is not 
of $\W_{\kappa}(\sg)$-modules.
However,
set 
\begin{align*}
\text{$G^p H_-^0(M(\lam)^*)=
D(\M(\gamma_{\bar{\lam}})/G^{-p}
\M(\gamma_{\bar{\lam}}))$
\quad ($p\geq -1$)
as $L(0)$-modules 
}
\end{align*}
via the isomorphism \eqref{eq:iso-dual-Verma-fact}.
Then,
\begin{align*}
 &H_-^0(M(\lam)^*)
=
G^{-1} H_-^0(M(\lam)^*)\supset 
G^0  H_-^0(M(\lam)^*)\supset \dots
\supset 
 G^p 
H_-^0(M(\lam)^*)\supset \dots, \\
& \bigcap_p G^p H_-^0(M(\lam)^*)=\{0\},
\end{align*}
and we see from the proof of \cite[Theorem 6,8]{A}
 that
\begin{align}
&G^p\W_{\kappa}(\sg)\cdot G^q H_-^0(M(\lam)^*)
\subset G^{p+q}H_-^0(M(\lam)^*),
\end{align}
and that 
\begin{align}
\gr^G H_-^0(M(\lam)^*)\cong \gr^G D(\M(\gamma_{\bar{\lam}}))
\end{align}
as modules over the polynomial ring generated by
$\bar{W}_i(n)$  $(i\in \sI, n>0)$.
Here,
of course, 
$\gr^G H_-^0(M(\lam)^*)
=\bigoplus_{p}G^p H_-^0(M(\lam)^*)/G^{p+1}H_-^0(M(\lam)^*)
$.
This shows, by \eqref{eq:singular-vector-of-dual-of-Verma},
that
\begin{align*}
\{v\in \gr^G H_-^0(M(\lam)^*);
\bar{W}_i(n)v=0\quad (i\in \sI, n>0)\}
=
\C \overline{|\lam}\dual{\ket}.
\end{align*}
Here,
 $\overline{|\lam}\dual{\ket}$
is the image of $|\lam\dual{\ket}$.
Hence
$\HW(H_-^0(M(\lam)^*))= \C |\lam \dual{\ket}$.
(1) is proved.
(2) can be similarly proved by using \cite[Theorem 6.8 (1)]{A}.
\end{proof}
\subsection{The generic Verma modules}
Let $\kappa\in \C^*$.
\begin{Th}\label{Th:generc--Verma}
Suppose that
$\lam\in \dual{\h}_{\kappa}$ 
is antidominant,
i.e,
$\bra \lam+\rho,\alpha\che\ket\not\in
\{1,2,\dots,\}
$ for all $\alpha\in \prroots$.
Then,
$\M(\gamma_{\bar{\lam}})=\why(\gamma_{\bar{\lam}})$.
\end{Th}
\begin{proof}
We have:
$M(\lam)=M(\lam)^*=L(\lam)$
for an antidominant $\lam$.
But then,
by Theorem \ref{Th:image_of_Verma_Module} and 
Theorem \ref{Th:dual},
$\M(\gamma_{\bar{\lam}})=D(\M(\gamma_{-w_0(\bar{\lam})}))$.
But by Lemma \ref{Lem:criterion-simple},
this only happens when
$\M(\gamma_{\bar{\lam}})
=\why(\gamma_{\bar{\lam}})$.
\end{proof}
\begin{Rem}
 Let $\sg={\mathfrak{sl}_2}(\C)$.
Then,
one can apply
Theorem \ref{Th:generc--Verma}
to give  yet another proof of
\cite[Proposition 8.2 $(b)$]{KR}.
\end{Rem}
\subsection{The functor $H_{\pm}^0(?)$}
\begin{Def}$ $
 
\begin{enumerate}
  \item A infinitesimal character
$\gamma_{\bar{\lam}}$ ($\bar{\lam}\in \dual{\sh}$)
is called {\em non-degenerate}
if 
$\bra \bar{\lam},{\alpha}\che\ket \not
\in \Z$ for all ${\alpha}\in \sroots$.
  \item A weight ${\Lam}\in {\h}$
is called {\em non-degenerate}
if $\bra {\Lam},\bar{\alpha}\che\ket \not
\in \Z$ for all $\bar{\alpha}\in \sroots$.
 \end{enumerate}
\end{Def}
\begin{Rem}$ $

 \begin{enumerate}
  \item $\Lam\in \dual{\h}$
is non-degenerate if and only if
$R^{\Lam}\cap \sroots=\emptyset$.
Therefore,
$\Lam$ is non-degenerate if and only if
$w\circ \Lam$
is non-egenerate for all $w\in W^{\Lam}$.
\item 
The integral Weyl group $W^{\Lam}$
can be an infinite group even when $\Lam\in \dual{\h}$
is non-degenerate.
Indeed,
$W^{\Lam}\cong W$
for a  principal admissible weight $\Lam$,
see \cite{KW1, KW2}.
 \end{enumerate}
\end{Rem}
\begin{Th}[{\cite[Theorem 8.3]{A}}]\label{Th:vanishing}
Let $\kappa\in \C^*$.

\begin {enumerate}
 \item Suppose $\Lam\in \dual{\h}_{\kappa}$ is non-degenerate.
Then,
$H_-^i(V)=0$
$(i\ne 0)$
for all $V\in \BGG_{\kappa}^{[\Lam]}$.
 \item Suppose $\Lam\in \dual{\h}_{\kappa}$ satisfies the following
condition.
\begin{align}\label{eq:cond+}
R^{\Lam}\cap \prroots\cap t_{\srho\che}(\nrroots)
=\emptyset. 
\end{align}
Then,
$H_+^i(V)=0$
$(i\ne 0)$
for all $V\in \BGG_{\kappa}^{[\Lam]}$.
\end{enumerate}
\end{Th}
\begin{Rem}\label{Rem:+and-}
$ $

 \begin{enumerate}
  \item 
The condition \eqref{eq:cond+}
is equivalent to
\begin{align}
\bra \Lam,\alpha\che\ket\not\in \Z
\quad \text{for all 
$\alpha\in \{-\bar{\alpha}+n\delta;
\bar{\alpha}\in \sproots,1\leq n\leq \height\bar{\alpha}\}$}.
\end{align}
\item If $\Lam\in \dual{\h}$
satisfies \eqref{eq:cond+},
then
$t_{-\srho\che}\circ \Lam$
is non-degenerate.
 \end{enumerate}
\end{Rem}
\begin{Co}\label{Co:exact}
Let $\kappa\in \C^*$.

\begin{enumerate}
 \item Suppose $\Lam\in \dual{\h}_{\kappa}$ is non-degenerate.
Then,
the correspondence $V\fmap H_{-}^0(V)$ 
defines an exact functor from 
$\BGG_{\kappa}^{[\Lam]}$
to $\Wcat$.
 \item Suppose $\Lam\in \dual{\h}_{\kappa}$ satisfies
\eqref{eq:cond+}.
Then,
the correspondence $V\fmap H_{+}^0(V)$ 
defines an exact functor from 
$\BGG_{\kappa}^{[\Lam]}$
to $\Wcat$.
\end{enumerate} 
\end{Co}
\subsection{The images of simple modules ($(-)$-case)}
\begin{Th}\label{Th:irr-}
Suppose
that
$\lam\in \dual{\h}$
is non-degenerate
and non-critical.
Then,
$H_-^0(L(\lam))\cong
\why(\gamma_{\bar{\lam}})$.
\end{Th}
\begin{proof}
By Corollary \ref{Co:exact},
Theorem \ref{Th:image_of_Verma_Module}
and Theorem \ref{Th:dual},
the exact sequences
$M(\lam)\rightarrow L(\lam)\rightarrow 0$
and $0\rightarrow L(\lam)
\rightarrow M(\lam)^*$
in $\BGG_{\kappa}^{[\lam]}$
give rise to
the exact sequences
\begin{align}
&\M(\gamma_{\bar{\lam}})\rightarrow
H_-^0(L(\lam))\rightarrow 0,
\label{eq:surj}
\\
&0
\rightarrow H_-^0(L(\lam))
\rightarrow D(\M(\gamma_{-w_0(\bar{\lam})})).
\label{eq:inj}
\end{align}
Thus,
by Lemma \ref{Lem:criterion-simple},
$H_-^0(L(\lam))$
must be  either $\{0\}$ or isomorphic to $\why(\gamma_{\bar{\lam}})$.
But
\cite[Proposition 7.6]{A}
shows that
$H_{-}^0(L(\lam))\ne \{0\}$.
\end{proof}
Theorem \ref{Th:irr-}
in particular proves the irreducibility conjecture
of Frenkel, Kac and Wakimoto,
see
\cite[Conjecture 3.4$_- (a)$]{FKW}.
One can also apply Theorem \ref{Th:vanishing} and 
Theorem \ref{Th:irr-} to
non-principal admissible weights to
prove the conjecture \cite[Proposition 3.4 $(c)$]{FKW}.

\subsection{The characters}
Let
\begin{align}
& \fdomain_{\pm}^{\kappa}=
\{\Lam \in \dual{\h}_{\kappa};
\bra \Lam+\rho, \alpha\che\ket
\not \in \{\mp 1,,\mp 2,\dots,\}
\text{ for all }
\alpha\in \prroots\},
\\
& \fdomain_{\pm, \nondeg}^{\kappa
}=\{
\Lam \in\fdomain_{\pm}^{\kappa};
\Lam \text{ is non-degenerate} \}
\subset \fdomain_{\pm}^{\kappa}.
\end{align}

For $\Lam\in \dual{\h}$,
let
$\geq_{\Lam}$
denote the Bruhat ordering in $W^{\Lam}$
and let $\ell_{\Lam}$
denote the length function on $W^{\Lam}$.
Let
$P_{w,y}^{\Lam}$  $(w,y \in W^{\Lam})$
be the corresponding Kazhdan-Lusztig polynomial
and $Q_{w,y}^{\Lam}$
the inverse 
Kazhdan-Lusztig polynomial.
Let
$W^{\Lam}_0=\bra s_{\alpha};
\bra \Lam+\rho,\alpha\che \ket =0\ket \subset W^{\Lam}$.

The following theorem 
describes the normalized characters
of all irreducible representations of $\W_{\kappa}(\sg)$
of non-generate highest weight $\gamma_{\bar{\lam}}$.
\begin{Th}\label{Th:ch-formula}
Let $\kappa \in \C^*$.

 \begin{enumerate}
  \item Let $\Lam \in \fdomain_{+, \nondeg}^{\kappa}$.
Then,
for any $w\in W^{\Lam}$
which is the longest element of $w W^{\Lam}_0$
we have
\begin{align*}
 \ch(\why(\gamma_{\overline{w\circ \Lam}}))
=\frac{1}{\eta(\tau)^{\rank \sg}}
\sum_{W^{\Lam} \ni
y\geq_{\Lam} w }(-1)^{\ell_{\Lam}(y)-\ell_{\Lam}(w)}
Q_{w,y}^{\Lam}(1)
{q^{\frac{|y(\Lam+\rho)|^2}{2\kappa}}}.
\end{align*}
  \item Let $\Lam \in \fdomain_{-, \nondeg}^{\kappa}$.
Then,
for any $w\in W^{\Lam}$
which is the shortest  element of $w W^{\Lam}_0$
we have
\begin{align*}
 \ch(\why(\gamma_{\overline{w\circ \Lam}}))
=\frac{1}{\eta(\tau)^{\rank \sg}}
\sum_{W^{\Lam} \ni
y\leq_{\Lam} w }(-1)^{\ell_{\Lam}(y)-\ell_{\Lam}(w)}
P_{w,y}^{\Lam}(1)
{q^{\frac{|y(\Lam+\rho)|^2}{2\kappa}}}.
\end{align*}
 \end{enumerate}
\end{Th}
\begin{proof}
By Corollary \ref{Co:exact} and Theorem \ref{Th:irr-},
Theorem  follows  directly from the character formula
\cite[Theorem 1.1]{KT}
in $\BGG_{\kappa}$.
\end{proof}

\subsection{The images of simple modules ($(+)$-case)}
\begin{Th}\label{Th:simple+}
Suppose
that
$\lam\in \dual{\h}$
is non-critical and
satisfies the condition \eqref{eq:cond+}.
Then,
$H_+^0(L(\lam))
\cong \why(\gamma_{\overline{t_{-\srho\che}\circ \lam}})$.
\end{Th}
\begin{proof}
By 
 Theorem \ref{Th:dual} (2)
and Corollary \ref{Co:exact} (2),
the condition on  $\lam$
implies that 
$H_+^0(L(\lam))$
is a submodule of
$D(\M(\gamma_{-w_0(\overline{
t_{-\srho\che}\circ \lam})}))$.
Therefore,
it is sufficient 
to show that
$
\ch H_+^0(L(\lam))=\ch
\why(\gamma_{\overline{t_{-\srho\che}\circ \lam}})
$.
But this follows from
\cite[Remark 3.5]{A},
 Remark \ref{Rem:+and-} (2),
Corollary \ref{Co:exact}
and Theorem \ref{Th:irr-}.
\end{proof}

\subsection{The simple vertex operator algebra}
Let 
\begin{align*}
 \text{$\vac=\overline{t_{-\srho\che}\circ
(\kappa-h\che)\Lam_0}=-\kappa \srho\che\in \dual{\sh}$.
}
\end{align*}

\begin{Pro}\label{Pro:simple-before}
There exist a unique surjection 
$\M(\gamma_{\vac})\twoheadrightarrow \W_{\kappa}(\sg)$
of $\W_{\kappa}(\sg)$-modules
that sends $|\gamma_{\vac}\ket$ to $|0\ket$.
\end{Pro}
\begin{proof}
By Theorem \ref{Th:image_of_Verma_Module},
 it is sufficient to show that
$\C |0\ket \cong \C |t_{-\srho\che}\circ (\kappa-h\che)\Lam_0\ket$
as $\A(\W_{\kappa}(\sg))$-modules.
The proof of this assertion is 
the same as the proof of
Lemma \ref{Lem:hw-Verma--} (2).
%
\end{proof}
Proposition \ref{Pro:simple-before}
implies the following result:
\begin{Th}
For any $\kappa\in \C^*$,
$\why(\gamma_{\vac})$ 
is the unique simple quotient of $\W_{\kappa}(\sg)$.
In particular,
$\why(\gamma_{\vac})$
carries a vertex operator algebra  structure.
\end{Th}
Of course the character of $\why(\gamma_{\vac})$
is expressed by
Theorem 
\ref{Th:ch-formula}
when $\vac\in \dual{\sh}$
is non-degenerate,
i.e,
when $\kappa(\srho\che,\alpha\che)\not\in \Z$
for all $\alpha\in \sroots$.
\subsection{The conjecture }
\begin{Conj}
$ $

\begin{enumerate}
\item For any $\kappa\in \C$ and any $V\in \Obj \BGG_{\kappa}$,
$H^i_-(V)=0$ ($i\ne 0$).
\item For any non-critical weight $\lam$,
\begin{align*}
H^0_-(L(\lam))\cong \begin{cases}
\why( \gamma_{\bar{\lam}})&\text{if 
$\bra \lam+\rho,\bar{\alpha}\che\ket\not\in \{1,2,\dots,\}$
for all $\bar{\alpha}\in \sproots$},\\
0&\text{otherwise.}
			 \end{cases}
\end{align*}
\end{enumerate}
\end{Conj}
%

\end{document}